\input amstex
\documentstyle{amsppt}
\pageheight{194mm}
\pagewidth{133mm}
\magnification\magstep1


\def\nologo{\let\logo@\empty}

\def\Ad{\operatorname{Ad}}

\def\Aut{\operatorname{Aut}}
\def\Aut{\operatorname{Aut}}

\def\BS{\operatorname{BS}}

\def\Cont{\operatorname{Cont}}

\def\deg{\operatorname{deg}}
\def\End{\operatorname{End}}

\def\cExt{\operatorname{\Cal Ext}}

\def\gr{\operatorname{gr}}
\def\Hom{\operatorname{Hom}}

\def\Im{\operatorname{Im}}

\def\Ker{\operatorname{Ker}}
\def\Lie{\operatorname{Lie}}

\def\LMH{\operatorname{LMH}}

\def\Map{\;\operatorname{Map}}

\def\rank{\operatorname{rank}}
\def\Re{\operatorname{Re}}

\def\Res{\operatorname{Res}}

\def\SL{\operatorname{SL}}

\def\Spec{\operatorname{Spec}}

\def\spl{\operatorname{spl}}

\def\Sym{\operatorname{Sym}}
\def\toric{\operatorname{toric}}

\def\torus{\operatorname{torus}}

\def\ts{\tsize\sum}

\def\bC{\bold C}

\def\be{\bold e}

\def\bN{\bold N}
\def\bP{\bold P}
\def\bQ{\bold Q}

\def\bR{\bold R}
\def\bS{\bold S}

\def\bZ{\bold Z}

\def\cB{{\Cal B}}

\def\cE{{\Cal E}}
\def\cF{{\Cal F}}

\def\cH{{\Cal H}}
\def\cI{{\Cal I}}

\def\cK{{\Cal K}}
\def\cL{{\Cal L}}
\def\cM{{\Cal M}}
\def\cN{{\Cal N}}
\def\cO{{\Cal O}}

\def\cT{{\Cal T}}

\def\fg{{\frak g}}

\def\a{\alpha}
\def\b{\beta}

\def\g{\gamma}
\def\G{\Gamma}

\def\sig{\sigma}
\def\Sig{\Sigma}

\def\.{$.\;$}

\def\gp{{\text{\rm gp}}}
\def\loga{{\text{\rm log}}}
\def\mult{{\text{\rm mult}}}

\def\resp.{\text{resp}.\;}

\def\tra{\overset{\sim}\to{\to}}
\def\tla{\overset{\sim}\to{\leftarrow}}
\def\triv{{\text{\rm triv}}}

\def\val{{\text{\rm val}}}

\def\O^logten{\cO\^log\otimes}

\let\bs=\backslash

\let\la=\leftarrow

\let\lan=\langle
\let\lan=\langle

\let\ox=\otimes

\let\ran=\rangle

\let\sub=\subset

\let\t=\tilde

\let\x=\times

\def\Dc{\check{D}}
\def\Ec{\check{E}}


\topmatter

\title
Studies of closed/open mirror symmetry for quintic threefolds through log mixed Hodge theory\\
\endtitle

\dedicatory
Dedicated to Kazuhiro Konno and dedicated  to James D.\ Lewis on his sixtieth birthday
\enddedicatory

\author
Sampei Usui\footnote{\text{Partially supported by JSPS. KAKENHI  (B) No. 23340008.}}
\endauthor

\date
\enddate

\address
\newline
{\rm Sampei USUI
\newline
Graduate School of Science
\newline
Osaka University
\newline
Toyonaka, Osaka, 560-0043, Japan}
\newline
{\tt usui\@math.sci.osaka-u.ac.jp}
\endaddress

\footnote"{}"{2010 {\it Mathematics Subject Classification}.
Primary 14C30; Secondary 14D07, 32G20.}

\abstract
We correct the definitions and descriptions  of the integral structures in [U14].
The previous flat basis in [ibid] is characterized by the Frobenius solutions and is integral in the first approximation by mean of the graded quotients of monodromy filtration, but it is not integral in the strict sense.
In this article, we use $\hat\Gamma$-integral structure of Iritani in [I11] for A-model.
Using this precise version, we study open mirror symmetry for quintic threefolds through log mixed Hodge theory, especially the recent result on N\'eron models for admissible normal functions with non-torsion extensions in the joint work [KNU14] with K. Kato and C. Nakayama.
We understand asymptotic conditions as values in the fiber over a base point on the boundary of $S^{\log}$.
\endabstract

\endtopmatter

\NoRunningHeads

\document

0. Introduction

1. Log mixed Hodge theory

\hskip10pt 
1.1. Category $\cB(\log)$

\hskip10pt 
1.2. Ringed space $(S^{\log}, \cO_S^{\log})$

\hskip10pt 
1.3. Toric variety

\hskip10pt 
1.4. Local systems on $(S^{\log}, \cO_S^{\log})$

\hskip10pt 
1.5. Graded polarized log mixed Hodge structure

\hskip10pt 
1.6. Nilpotent orbit

\hskip10pt 
1.7. Moduli of log mixed Hodge structures of specified type

\hskip10pt 
1.8. N\'eron model for admissible normal function

2. Quintic threefolds

\hskip10pt 
2.1. Quintic mirror family

\hskip10pt 
2.2. Picard--Fuchs equation on B-model of mirror $V^\circ$

\hskip10pt 
2.3. A-model of quintic $V$

\hskip10pt 
2.4. Integral structure

\hskip10pt 
2.5. Correspondence table

\hskip10pt 
2.6. Proofs of results in 2.5

3. Proof of (5) in Introduction

\hskip10pt 
3.1. Proof of (5.1)--(5.3) in Introduction over log disc $S$

\hskip10pt 
3.2. Proof of (5.1)--(5.2) in Introduction over log point $p_0$

\hskip10pt 
3.3. Discussions on geometries and local systems

\bigskip


{\bf 0. Introduction}
\medskip

In a series of joint works with K. Kato and C. Nakayama, we are constructing a fundamental diagram which consists of various kind of partial compactifications of classifying space of mixed Hodge structures and their relations.
We try to understand mirror symmetry in this framework of the fundamental diagram.
In this paper, we first complete the insufficient results in the previous paper of Usui [U14, 3.5--3.6] (see Remark in 2.6 below), and then study open mirror symmetry for quintic threefolds through log mixed Hodge theory, especially the fine moduli of log Hodge structures and N\'eron models over it.
\medskip

{\it Fundamental Diagram}

For a classifying space $D$ of Hodge structures of specified type, we have
$$
\CD
@.@.D_{\SL(2),\val}@.@>>>D_{\BS,\val}\\
@.@.@VVV@.@VVV\\
\Gamma\bs D_{\Sig,\val}@<<<D_{\Sig,\val}^\sharp@>>>D_{\SL(2)}@.@.D_{\BS}\\
@VVV@VVV\\
\Gamma\bs D_{\Sig}@<<<D_{\Sig}^\sharp
\endCD
$$
in pure case: [KU99], [KU02], [KU09].
For mixed case, we should extend to an amplified diagram: 
[KNU08], [KNU09], [KNU11], [KNU13], continuing.

\medskip

{\it Mirror symmetry for quintic threefolds}
\medskip

Let $V$ be a quintic threefold in $\bP^4$ and $V_\psi^\circ$ be its mirror family (cf.\ [CoK99, Sect.\ 4.2]).
For simplicity, we denote the family $V_\psi^\circ$ simply by $V^\circ$ if there would be no confusions.

Mirror symmetry for the A-model of quintic threefold $V$ 
and the B-model of its mirror $V^\circ$ was predicted by Candelas--de la Ossa--Green--Parks  
in the famous paper [CDGP91].
We recall two styles of the theorem (1) and (2) below.
Every statement in the present paper is near the large radius point $q_0$ of the complexified K\"ahler moduli $\cK\cM(V)$ and the maximally unipotent monodromy point $p_0$ of the complex moduli $\cM(V^\circ)$.

\medskip

Let $t:=y_1/y_0$, $u:=t/2\pi i$ be the canonical parameters and $q:=e^t=e^{2\pi iu}$ be the canonical coordinate for B-model in 2.2 below and the respective ones for A-model in 2.3 below.

The following theorem is due to Givental [G96] and Lian--Liu--Yau [LLuY97].

\medskip

\noindent
{\bf (1)} ({\it Potential}).
{\it The potentials of the two models coincide: 
$\Phi^V_{\text{\rm GW}}(t)=\Phi^{V^\circ}_{\text{\rm GM}}(t)$.}
\medskip

Morrison [M97] formulated the following style (2) and proved the theorem except integral structure.
Iritani [I11] defined $\hat\Gamma$-integral structure for A-model and proved the theorem completely for wider objects.
\medskip

\noindent
{\bf (2)} ({\it Variation of Hodge structure}).
{\it The isomorphism $(q_0\in\overline{\cK\cM}(V))\tla(p_0\in\overline\cM(V^\circ))$ of neighborhoods of the compactifications, by the canonical coordinates $q=\exp(2\pi iu)$, lifts to an isomorphism, over the punctured neighborhoods $\cK\cM(V)\tla\cM(V^\circ)$, of polarized $\bZ$-variations of Hodge structure with a specified section  
$$
(\cH^V,S, \nabla^{\text{even}},\cH_\bZ^V,F;T^3)
\tla(\cH^{V^\circ},Q, \nabla^{\text{GM}},\cH_\bZ^{V^\circ},F;\t\Omega).
$$
}
\medskip

Our (3) below is equivalent to (1) and (2) by a log version [KU09, 2.5.14] of the nilpotent orbit theorem of Schmid [S73] (this part of [U14] is valid).
\medskip

\noindent
{\bf (3)} ({\it Log Hodge structure, Log period map}).
{\it The isomorphism $(q_0\in\overline{\cK\cM}(V))\tla(p_0\in\overline\cM(V^\circ))$ of neighborhoods of the compactifications uniquely lifts to an isomorphism of B-model  log variation of polarized Hodge structure with a specified section $\t\Omega$ for $V^\circ$ and A-model log variation of polarized Hodge structure with a specified section $T^3$ for $V$, whose restriction over the punctured $\cK\cM(V)\tla\cM(V^\circ)$ coincides with the isomorphism of variations of polarized Hodge structure with specified sections in (2).

This rephrases as follows.
Let $\sigma$ be the common monodromy cone, transformed by a level structure into $\End$ of a reference fiber of the local system, for the A-model and for the B-model.
Then, we have a commutative diagram of horizontal log period maps
$$
\CD
(q_0\in\overline{\cK\cM}(V))\ \tla\ (p_0\in\overline\cM(V^\circ))\\
\searrow\qquad\qquad\swarrow\\
([\sig,\exp(\sig_\bC)F_0]\in\Gamma(\sigma)^{\gp}\bs D_\sig)
\endCD
$$
with extensions of specified sections in (2), where $(\sig,\exp(\sig_\bC)F_0)$ is the nilpotent orbit, regarded as a boundary point, and $\Gamma(\sigma)^{\gp}\bs D_\sig$ is the fine moduli of log Hodge structures of specified type which will be explained in Section 1 below.}
\bigskip

{\it Open mirror symmetry for quintic threefolds}
\medskip

The following theorem is due to Walcher [W07] and Morrison--Walcher [MW09].
\medskip

\noindent
{\bf (4)} ({\it Inhomogenous solutions}).
{\it Let $\cL$ be the Picard--Fuchs differential operator for quintic mirror family (cf.\ 2.2 below). 
Let
$$
\cT_A=\frac{u}{2}\pm\Big(\frac{1}{4}+\frac{1}{2\pi^2}\sum_{d\;\text{odd}}n_dq^{d/2}\Big)
$$
be the A-model domainwall tension in [MW09], where the $n_{d}$ are open Gromov--Witten invariants, and 
$$
\cT_B=\int_{C_-}^{C_+}\Omega
$$
be the B-model domainwall tension, where $C_{\pm}\sub V^{\circ}$ are the disjoint smooth curves coming from the two conics in $\{x_1+x_2=x_3+x_4=0\}\cap V_\psi\subset\bP^4(\bC)$ [ibid].

Then
$$
\cL(y_0(z)\cT_A(z))=\cL(\cT_B(z))\Big(=\frac{15}{16\pi^2}\sqrt{z}\Big)\quad 
\Big(z=\frac{1}{(5\psi)^5}\Big).
$$
}
\medskip

Concerning this, we have the following observations.
\medskip

\noindent
{\bf (5)} ({\it Log mixed Hodge structure, Log normal function}).
We describe for B-model. 
The same holds for A-model by (1)--(4) and the correspondence table in 2.5 below.

Put $\cH:=\cH^{V^\circ}$ and $\cT:=\cT_B$.
We use $e^0\in I^{0,0}$, $e^1\in I^{1,1}$ which are a part of a basis of 
$\cH_{\cO}$ respecting the Deligne decomposition at $p_0$ (see 2.5 (3B)) and a part of flat basis $s^0=e^0$, $s^1=e^1-ue^0$ of $\cH_\bZ$ (see 2.5 (5B)).
To make the local monodromy of $\cT$ unipotent, we take a double cover $z^{1/2}\mapsto z$.
Let $L_\bQ$ be the translated local system comparing to the trivial extension $\bQ(-2)\oplus \cH_\bQ$ by $(0,-(\cT/y_0)s^0)$ in $\cE xt^1(\bQ(-2), \cH_\bQ)$.
Let  $J_{L_\bQ}$ be the N\'eron model on a neighborhood $S$ of $p_0$ in the $z^{1/2}$-plane which lies over $L_\bQ$ in [KNU14] (there is a difference of Tate twist).
Then, $J_{L_\bQ}=\cE xt^1_{\text{LMH}/S}(\bZ(-2), \cH)$ (extension group of log mixed Hodge structures over $S$) in the present case (cf.\ [KNU13, III, Corollary 6.1.6], and 1.8 below), and we have the following (5.1)--(5.3).
\medskip

\noindent
{\bf (5.1)} 
{\it The normalized tension $\cT/y_0$ is understood as a multi-valued truncated normal function by $(\cT/y_0)s^0$.
Then it lifts and extends uniquely to a single-valued log normal function $S\to J_{L_\bQ}$ 
so that the corresponding exact sequence $0\to\cH\to\tilde\cH\to\bZ(-2)\to0$ of log mixed Hodge structures over $S$ is given by the liftings $1_\bZ$ and $1_F$ in $\tilde\cH$ of $(2\pi i)^{-2}\cdot1\in\bZ(-2)$ respecting the lattice and the Hodge filtration, respectively, which are defined as follows:
$1_\bZ:=((2\pi i)^{-2}\cdot1,\,-(\cT/y_0)s^0)$ with $(\cT/y_0)s^0\in \cH_{\cO^{\log}}=(\gr^W_3)_{\cO^{\log}}$, and 
$1_F-1_\bZ:=(\delta(\cT/y_0))e^1-(\cT/y_0)e^0$, where $\delta:=2\pi iqd/dq$.}
\medskip

\noindent
{\bf (5.2)} 
{\it A splitting of the weight filtration $W$ of the local system $\tilde\cH_\bZ$, i.e., a splitting 
compatible with the monodromy of the local system $\tilde\cH_\bZ$, is given by
$1_\bZ^{\text{spl}}=1_\bZ+s^1/2$, and the log normal function over it is given by
$1_F^{\text{spl}}-1_\bZ^{\text{spl}}=(\delta(\cT/y_0))e^1-(\cT/y_0)e^0$, where $\delta$ is as in (5.1).}
\medskip

\noindent
{\bf (5.3)} 
{\it (4) says that the inverse of the normal function in (5.1) from its 
image is given by $16\pi^2/15$ times the Picard--Fuchs differential operator $\cL$.}
\medskip

(5) is proved in Section 3, and after these proofs some geometric backgrounds are discussed in Section 3.3.
\medskip

The organization of this paper is as follows.
Section 1 is a summary of log mixed Hodge theory mainly from [KU09], [KNU13] and [KNU14], which is used to study mirror symmetry in later sections and also is expected to work as a brief guide of this theory.
In Section 2, after preparations including $\hat\Gamma$-integral structure in [I09] and  [I11], we give a correspondence table of closed mirror symmetry for quintic threefolds and their mirrors, which is the precision and the expansion of our previous paper [U14, 3].
In Section 3, we prove log mixed Hodge theoretic interpretation (5).
We also give some discussions on the related geometries and local systems in Section 3.3.

\medskip

{{\it Acknowledgments.}
The author thanks Kazuya Kato and Chikara Nakayama for series of joint works on log Hodge theory, from which he learns a lot and enjoys  exciting studies.
He thanks Hiroshi Iritani for pointing out insufficient parts in the previous paper [U14].
He also thanks Yukiko Konishi and Satoshi Minabe, together with Iritani, for stimulating seminars on the present topic.
The author thanks the referees for careful readings and useful comments.

\medskip

{\it Notation.}

Fix $\Lambda:=(H_0,W,(\lan\;,\;\ran_w)_w,(h^{p,q})_{p,q})$, where 

$H_0$ is a free $\bZ$-module of finite rank, 

$W$ is an increasing filtration on $H_{0,\bQ}:=\bQ\otimes H_0$, 

$\lan\;,\;\ran_w$ is a non-degenerate $(-1)^w$-symmetric bilinear form on $\gr^W_w$, 

$(h^{p,q})_{p,q}$ is a set of Hodge numbers.

\medskip

$D$ : the classifying space of graded polarized mixed Hodge structures for the data $\Lambda$, consisting of all Hodge filtrations.

$\Dc$ : the \lq\lq compact dual" of $D$.

$G_A:=\Aut(H_{0,A},W,(\lan\;,\;\ran_w)_w)$, where $H_{0,A}:=A\otimes H_0$  
$(A=\bZ,\bQ,\bR,\bC)$,

$\fg_A:=\Lie G_{A}=\End(H_{0,A},W,(\lan\;,\;\ran_w)_w)$ 
$(A=\bQ,\bR,\bC)$.
\medskip

We treat Tate twists case by case in  this paper.
We hope to treat them consistently in a future.

\bigskip

{\bf 1. Log mixed Hodge theory}
\medskip

This section is a summary of log mixed Hodge theory from [KU09], [KNU13] and [KNU14].
We write a general form of these results as a brief guide for future use.
Section 1.8 is adapted for the use in Section 3.
The corresponding results in [KNU13] and [KNU14] are written in more general settings.
\medskip

{\bf 1.1. Category $\cB(\log)$}
\medskip

Let $S$ be a subset  of an analytic space $Z$.
The {\it strong topology of $S$ in $Z$} is the strongest one among those topologies on $S$ in which, for any analytic space $A$ and any morphism $f:A\to Z$ with $f(A)\sub S$ as sets, $f:A\to S$ is continuous.
$S$ is regarded as a local ringed space by the pullback sheaf of $\cO_Z$.

Let $\cB$ be the category of local ringed spaces $S$ over $\bC$ which have an open covering $(U_\lambda)_\lambda$ satisfying the following condition:
For each $\lambda$, there exist an analytic space $Z_\lambda$, and a subset $S_\lambda$ of $Z_\lambda$ such that, as local ringed space over $\bC$, $U_\lambda$ is isomorphic to an open subset of $S_\lambda$ which is endowed with the strong topology in $Z_\lambda$ and the inverse image of $\cO_{Z_\lambda}$.

A {\it log structure} on a local ringed space $S$ is a sheaf of monoids $M$ on $S$ together with a homomorphisim $\alpha: M\to\cO_S$ such that $\alpha^{-1}\cO_S^\times\tra\cO_S^\times$.
{\it fs} log structure means, locally on the underlying space, the log structure has a chart 
which is finitely generated, integral and saturated.

Let $\cB(\log)$ be the category of objects of $\cB$ endowed with an fs log structure.

A log analytic space is called {\it log smooth} if, locally, it is isomorphic to an open set of a toric variety endowed with the canonical log structure.
A {\it log manifold} is a log local ringed space over $\bC$ which has an open 
covering $(U_\lambda)_\lambda$ satisfying the following condition:
For each $\lambda$, there exist a log smooth fs log analytic space $Z_\lambda$, 
a finite subset $I_\lambda$ of global log differential $1$-forms $\Gamma(Z_\lambda,\omega_{Z_\lambda}^1)$, and an isomorphism of log local ringed spaces over $\bC$ between $U_\lambda$ and an open subset of $S_\lambda:=\{z\in Z_\lambda\;|\;\text{the image of $I_\lambda$ in the stalk $\omega_z^1$ is zero}\}$ in the strong topology in $Z_\lambda$.
\bigskip

{\bf 1.2. Ringed space $(S^{\log}, \cO_S^{\log})$}
\medskip

The ringed space $(S^{\log}, \cO_S^{\log})$ was defined for fs log schemes by K.\ Kato and C.\ Nakayama in [KN99].
It was generalized for the category $\cB(\log)$ in [KU99].

Let $S\in\cB(\log)$.
As a set define

$S^{\log}:= \{(s,h) \,|\, \text{$s\in S$, $h: M_s^{\gp}\to\bS^1$ homomorphism
s.t.\ $h(u)=u/|u|$ if $u\in \cO_{S,s}^\times$}\}.$

\noindent
Endow $S^{\log}$ with the weakest topology such that the following two maps are continuous.

(1) $\tau: S^{\log}\to S, (s,h) \mapsto s$.

(2) For any open set $U\sub S$ and any $f\in \Gamma(U, M^{\gp})$, $\tau^{-1}(U)\to\bS^1$, $(s,h)\mapsto h(f_s)$.
\medskip

Then, $\tau$ is proper and surjective with fiber $\tau^{-1}(s)=(\bS^1)^{r(s)}$, where 
$r(s)$ is the rank of $(M^{\gp}/\cO_S^\times)_s$ which varies with $s\in S$.

Define a sheaf $\cL$ on $S^{\log}$ as the fiber product:
$$
\CD
\cL@>\exp>>\tau^{-1}(M^{\gp}) @.\quad\ni @.\quad(\text{$f$ at $(s,h)$})\\
@VVV@VVV@.\quad@VVV@.\\
\Cont(\ast, i\bR)@>\exp>>\Cont(\ast, \bS^1) @.\quad\ni @. \quad\,\; h(f).
\endCD
$$
Let $\iota: \tau^{-1}(\cO_S)\to\cL$ be a morphism  induced from
$$
\CD
f @. \quad\in@.\quad \tau^{-1}(\cO_S)@>\exp>> \tau^{-1}(\cO_S^{\times})\sub \tau^{-1}(M^{\gp})\\
@VVV@.\quad@VVV@VVV\\
(f-\bar f)/2@.\quad \in@.\quad\;\;\Cont(\ast, i\bR)@>\exp>>\;\Cont(\ast, \bS^1).
\endCD
$$
Define
$$
\cO_S^{\log}:=\frac{\tau^{-1}(\cO_S)\otimes\Sym_\bZ(\cL)}{(f\otimes1-1\otimes\iota(f)\,|\, f\in\tau^{-1}(\cO_S))}.
$$
Thus $\tau: (S^{\log}, \cO_S^{\log})\to (S, \cO_S)$ is a morphism of ringed spaces over $\bC$.
For $s\in S$ and $t\in S^{\log}$ lying over $s$, let $t_j\in\cL_t$ $(1\le j\le r(s))$ 
be elements such that their images in $(M^{\gp}/\cO_S^{\times})_s$ of $\exp(t_j)$ form a basis.
Then, $\cO_{S,t}^{\log}=\cO_{S,s}[t_j\;(1\le j\le r(s))]$ is a polynomial ring.

\bigskip

{\bf 1.3. Toric variety}
\medskip

Toric varieties offer typical examples of $S^{\log}$ and also they are building blocks of fine moduli spaces of log mixed Hodge structures.

Let $\sigma\sub\fg_\bR$ be a  {\it nilpotent cone}, i.e., a sharp cone, $\sigma\cap(-\sigma)=\{0\}$, generated by a finite number of mutually commutative nilpotent elements.
Assume that the cone generators of $\sig$ can be taken from $\fg_\bQ$.
Let $\Gamma$ be a subgroup of $G_\bZ$.
Define a monoid $\Gamma(\sigma):=\Gamma\cap\exp(\sigma)$ and the dual monoid $P(\sigma):=\G(\sig)^\vee=\Hom(\G(\sig),\bN)$.
Define a toric variety and a torus by
$$
\toric_\sig:=\Spec(\bC[P(\sig)])(\bC)=
\Hom(P(\sig),\bC^{\mult}) \supset\torus_\sig:=\Hom(P(\sig)^\gp,\bC^\times),
$$
where $\bC^\mult$ is $\bC$ regarded as a monoid by multiplication and 
$P(\sig)^\gp$ is the group generated by the monoid $P(\sig)$.
The exponential sequence $0\to\bZ\to\bC\to\bC^\times\to1$ induces the universal covering of the torus
$$
0\to\Hom(P(\sig)^\gp,\bZ)\to\Hom(P(\sig)^\gp,\bC)@>\be>>\Hom(P(\sig)^\gp,\bC^\times)\to1,
$$
where $\Hom(P(\sig)^\gp,\bZ)=\G(\sig)^{\gp}$ is considered as the fundamental group of $\torus_\sig$, and $\be(z\otimes\g):=e^{2\pi iz}\otimes\g$ ($z\in \bC$, $\g\in\G(\sig)^\gp=\Hom(P(\sigma)^{\gp},\bZ)$).

Fix the above cone $\sigma$.
For a face $\rho$ of $\sigma$, define $\tilde P(\rho):=\{l\in P(\sigma)^{\gp}\,|\,l(\rho)\ge0\}$.
Then we have an open covering
$$
\toric_\sig=\Spec(\bC[P(\sig)])(\bC)=\bigcup_{\rho\prec\sigma}\Spec(\bC[\tilde P(\rho)])(\bC).
$$

We now recall a stratification.
Fix the above cone $\sigma$ and let $\rho$ be a face of the cone $\sig$.
Then, we have a homomorphism $P(\rho)\la P(\sig)$ and hence a morphism $\toric_\rho\to\toric_\sig$.
The origin $0_\rho\in\toric_\rho$ is the monoid homomorphism $P(\sigma)\to\bC^{\mult}$ sending $1$ to $1$ and all the other elements of $P(\rho)$ to $0$, which is sent to a point of $\toric_\sig$ by the above morphism.
Then, as a set, we have a stratification into torus orbits
$$
\toric_\sig=\{\be(z)0_\rho\,|\,\rho\prec\sig,\, z\in\sig_\bC/(\rho_\bC+\log\G(\sig)^\gp)\}.
$$
Here $\be(c\log\gamma):=\be(c\otimes\gamma)=e^{2\pi ic}\otimes\gamma$ $(c\in \bC, \gamma\in \Gamma(\sigma)^{\gp})$.

For $S:=\toric_\sig$, the polar coordinate $\bR_{\ge0}\times\bS^1\to\bR_{\ge0}\bS^1=\bC$ induces
$\tau:S^{\log}\to S$ as 
$$
\align
\tau:S^{\log}&=\Hom(P(\sig),\bR_{\ge0}^\mult)\times\Hom(P(\sig),\bS^1)\\
&=\{(\be(iy)0_\rho,\be(x))\,|\,\rho\prec\sig,\, x\in\sig_\bR/(\rho_\bR+\log\G(\sig)^\gp),\,y\in\sig_\bR/\rho_\bR\}\\
&\to S=\Hom(P(\sig),\bC^{\mult}),\\
\tau(\be(ib)&0_\rho,\be(a))=\be(a+ib)0_\rho.
\endalign
$$
Since $0\to\rho_\bR/\log\G(\rho)^\gp\to\sig_\bR/\log\G(\sig)^\gp\to\sig_\bR/(\rho_\bR+\log\G(\sig)^\gp)\to0$ is exact, the fiber of $\tau$, as a set, is described as
$$
\tau^{-1}(\be(a+ib)0_\rho)
=\{(\be(ib)0_\rho,\be(a+x))\,|\,x\in\rho_\bR/\log\G(\rho)^\gp\}
\simeq(\bS^1)^r,
$$
where $r=r(\rho):=\rank\rho$ varies with $\rho\prec\sig$.

Let $H_\sig=(H_{\sig,\bZ},W,(\lan\;,\;\ran_w)_w)$ be the canonical local system endowed with the weight filtration and the polarizations on graded quotients on $S^{\log}$, which are given by the representation $\pi_1(S^{\log})=\G(\sig)^\gp\sub G_\bZ=\Aut(H_0,W,(\lan\;,\;\ran_w)_w)$.

\bigskip

{\bf 1.4. Local systems on $(S^{\log},\cO_S^{\log})$}
\medskip

We recall three results about local systems on $(S^{\log},\cO_S^{\log})\in\cB(\log)$ from [KU09, 2.3].

Let $L$ be a locally constant sheaf of abelian groups  on $S^\loga$.
For $s\in S$ and $t\in S^\loga$ lying over $s$, we call the action of $\pi_1(s^\loga)=\pi_1(\tau^{-1}(s))$ on $L_t$ the {\it local monodromy} of $L$ at $t$.
We say the local monodromy of $L$ is {\it unipotent} if the local monodromy of $L$ at $t$ is unipotent for any $t\in S^\loga$.

Let $s\in S$.
Let $(q_j)_{1\leq j\leq n}$ be a finite family of elements of
$M^\gp_{S,s}$ whose image in $(M^\gp_S/\cO^\x_S)_s$ is a
free basis, and let $(\g_j)_{1\leq j\leq n}$ be the dual
basis of $\pi_1(s^\loga)$, that is $[\gamma_j,q_k]=(2\pi i)\delta_{jk}$ where $[\;,\;]$ is the pairing given by 
$\pi_1(s^\loga)\simeq\Hom(M^\gp_s/\cO^\x_s,\bZ)$.

Let $L$ be a locally constant sheaf on $S^\loga$ of free $\bZ$-modules of finite rank.
Let $s\in S$ and $t\in \tau^{-1}(s)$, and
assume that the local monodromy of $L$ at $t$ is unipotent.
For a fixed $t$, we denote $L_0$ the constant sheaf on $S^\loga$ with fiber $L_t$.
Let $L_{0,\bQ}=\bQ\otimes_AL_0$, and let $N_j : L_{0,\bQ}\to L_{0,\bQ}$
be the endomorphism of constant sheaf which is induced by the logarithm of the monodromy action of $\g_j$ on the stalk $L_t$ of the locally constant sheaf $L$.
Lift $q_j$ in $\G(S, M^\gp_S)$ $($by replacing $S$ by an open
neighborhood of $s)$, and let
$$
\xi=\exp(\ts_{j=1}^n (2\pi i)^{-1}\log(q_j)\ox N_j):
\cO_S^\loga\ox_\bQ L_{0,\bQ}\tra\cO_S^\loga\ox_\bQ L_{0,\bQ}.
$$
Note that the operator $\xi$ depends on the choices of the branches of $\log(q_j)$ in $\cO_S^{\log}$ locally on $S^{\log}$, but that the subsheaf $\xi^{-1}(1\ox L_0)$ of $\cO_S^\loga \ox_\bZ L_0$ is independent of the choices and hence is defined globally on $S^\loga$.

The following proposition shows that the locally constant sheaf $L$ is embedded in $\cO_S^{\log}\otimes L_0$.

\proclaim{Proposition {\rm [KU09, Prop.\ 2.3.2]}}
Let the situation be as above.
If we replace $S$ by some open neighborhood of
$s$, we have an isomorphism of $\cO_S^\loga$-modules
$$
\nu: \cO_S^\loga \ox_A L \tra \cO_S^\loga \ox_A L_0
$$
satisfying the following condition $(1)$.

\noindent
{\rm(1)} The restriction of $\nu$ to $L=1\ox L$ induces an isomorphism
of locally constant sheaves $\nu  : L\tra \xi^{-1}(1\ox L_0)$.

If we take suitable branches $\log(q_j)_{t,0}$ in $\cO_{S,t}^{\log}$ of the germs  $\log(q_j)_t$ at $t$ $(1\leq j\leq n)$,
we can take an isomorphism $\nu$ which satisfies above $(1)$ and also the following $(2)$.

\noindent
{\rm(2)} The branch $\xi_{t,0}$ of the germ $\xi_t$, defined by the fixed branches $\log(q_j)_{t,0}$ of the germs $\log(q_j)_t$, satisfies
$\nu(1\otimes v)= \xi_{t,0}^{-1}(1\otimes v)$ for any $v\in L_t=L_0$.
\endproclaim

The following proposition yields a log Hodge theoretic understanding [ibid, 2.3.6] of the 
canonical extension of Deligne in [D70].

\proclaim{Proposition {\rm [KU09, Prop.\ 2.3.3]}}
Let $S\in\cB(\log)$ and let $L$ be a locally constant sheaf of finite dimensional
$\bC$-vector spaces on $S^\loga$.

{\rm(i)} If the local monodromy of $L$ is unipotent, the $\cO_S$-module
$\cM:=\tau_*(\cO^\loga_S\ox_\bC L)$
is locally free of finite rank, and we have an isomorphism
$\cO_S^\loga \ox_{\cO_S}\cM\tra\cO_S^\loga\ox_\bC L$.

{\rm(ii)} Conversely, assume that there are a locally free
$\cO_S$-module $\cM$ of finite rank on $S$ and an isomorphism
of $\cO_S^\loga$-modules $\cO_S^\loga\ox_{\cO_S}\cM \simeq \cO_S^\loga \ox_\bC L$.
Then the local monodromy of $L$ is unipotent and
$\cM\tra\tau_*(\cO_S^\loga\ox_\bC L)$.
\endproclaim

The following proposition describes the relation of log Gauss-Manin connection and monodromy logarithm over a log point.

\proclaim{Proposition {\rm [ibid, Prop.\ 2.3.4]}}
Let $S\in\cB(\log)$, let $L$ be a locally constant sheaf on $S^\loga$ of free $\bQ$-modules of finite rank.
Assume that the local monodromy of $L$ is unipotent.
\medskip

{\rm(i)} There exists a unique $\bQ$-homomorphism
$$
\cN : L \to (M_S^\gp/\cO_S^\x) \ox L
$$
satisfying the following condition $(1)$.
\medskip

\noindent
$(1)$ For any $s\in S$, any $t\in S^\loga$ lying over $s$, and any $\g\in \pi_1(s^\loga)$, if $h_\g: (M_S^\gp/\cO_S^\x)_s\to \bZ$
denotes the homomorphism corresponding to $\gamma$ by $\pi_1(s^\loga)\simeq\Hom(M^\gp_s/\cO^\x_s,\bZ)$, the
composition $L_t@>\cN>> (M_S^\gp/\cO^\x_S)_s\ox L_t@>h_\g>>
L_t$ coincides with the logarithm of the action of $\g$ on $L_t$.
\medskip

{\rm(ii)} Assume that $S$ is an fs log point $\{s\}$.
Let
$$
\cN' : L \to \omega^1_s\ox L
$$
be the composition of $\cN$ and the $\bQ$-linear map
$M_s^\gp/\cO_s^\x\ox L \to \omega^1_s\ox L, \;f\ox v\mapsto
(2\pi i)^{-1}d\log(f)\ox v$, and let $1\ox \cN': \cO_s^\loga \ox L
\to \omega^{1,\log}_s\ox L$ be the $\cO^\loga_s$-linear
homomorphism induced by $ \cN'$.
Let $\cM:=H^0(s^\loga,\, \cO^\loga_s\ox L)
=\tau_*(\cO^\loga_s\ox L)$.
Then the restriction $\cM\to \omega^1_s\ox_{\bC} \cM$ of
$d\ox 1_L : \cO^\loga_s\ox L \to \omega_s^{1,\loga}\ox L$
coincides with the restriction of $1\ox \cN'$ to $\cM$.
\endproclaim

$\cN$ in the above proposition is described as follows.
Assume $L=\xi^{-1}(1\otimes L_0)$ as in the first proposition.
Then $\cN(\xi^{-1}(1\otimes v)):=\sum_{j=1}^n q_j\otimes\xi^{-1}(1\otimes N_jv)$ for $v\in L_0$.

\bigskip

{\bf 1.5. Graded polarized log mixed Hodge structure}
\medskip

Let $S\in\cB(\log)$.
A {\it pre-graded polarized log mixed Hodge structure on $S$} is a tuple 
$H=(H_\bZ,W,(\lan\;,\;\ran_w)_w,H_\cO)$ consisting of a local system of $\bZ$-free modules $H_\bZ$ of finite rank on $S^{\log}$, an increasing filtration $W$ of $H_\bQ:=\bQ\otimes H_\bZ$, a non-degenerate $(-1)^w$-symmetric $\bQ$-bilinear form $\lan\;,\;\ran_w$ on $\gr^W_w$, a locally free $\cO_S$-module $H_\cO$ on $S$, 
a specified isomorphism $\cO_S^{\log}\otimes_\bZ H_\bZ\simeq \cO_S^{\log}\otimes_{\cO_S}H_\cO$ ({\it log Riemann-Hilbert correspondence}), and 
a specified decreasing filtration $FH_\cO$ of $H_\cO$ such that $F^pH_\cO$ and $H_\cO/F^pH_\cO$ are locally free.
Put $F^p:=\cO_S^{\log}\otimes_{\cO_S}F^pH_\cO$. 
Then $\tau_\ast F^p=F^pH_\cO$.
For each integer $w$, the orthogonality condition 
$\lan F^p(\gr^W_w),F^q(\gr^W_w)\ran_w=0$ $(p+q>w)$ is imposed.

A pre-graded polarized log mixed Hodge structure on $S$ is a {\it graded polarized log mixed Hodge structure on $S$} if its pullback to 
each $s\in S$ is a graded polarized log mixed Hodge structure on $s$ in the following sense.

Let $(H_\bZ,W,(\lan\;,\;\ran_w)_w,H_\cO)$ be a pre-graded polarized log mixed Hodge structure on a log point $s$.
It is a {\it graded polarized log mixed Hodge structure} if it satisfies the following three conditions.

(1) (Admissibility). 
For each logarithm $N$ of the local monodromy of the local system $(H_\bR, W,(\lan\;,\;\ran_w)_w)$, there exists a $W$-relative $N$-filtration $M(N,W)$.

(2) (Griffiths transversality). 
For any integer $p$, $\nabla F^p\sub \omega_s^{1,\log}\otimes F^{p-1}$ is satisfied, 
where $\omega_s^{1,\log}$ is the sheaf of $\cO_s^{\log}$-module of log differential 1-forms on $(s^{\log},\cO_s^{\log})$, and
$\nabla=d\otimes1_{H_\bZ}:\cO_s^{\log}\otimes H_\bZ\to\omega_s^{1,\log}\otimes H_\bZ$ is the log Gauss-Manin connection.

(3) (Positivity). 
For a point $t\in s^{\log}$ and a $\bC$-algebra homomorphism $a:\cO_{s,t}^{\log}\to\bC$, define a filtration $F(a):=\bC\otimes_{\cO_{s,t}^{\log}}F_t$ on $H_{\bC,t}$.
Then, $(H_{\bZ,t}(\gr^W_w),\lan\;,\;\ran_w,F(a))$ is a polarized Hodge structure  of weight $w$ in the usual sense if $a$ is sufficiently twisted, i.e., 
for $(q_j)_{1\le j\le n}\sub M_s$ inducing generators of $M_s/\cO_s^\times$, 
$|\exp(a(\log q_j))|\ll1$ for any $j$.

\bigskip

{\bf 1.6. Nilpotent orbit}
\medskip

Let $\sigma\sub\fg_\bR$ be a nilpotent cone (see 1.3).
A subset $Z\sub\Dc$ is {\it $\sig$-nilpotent orbit} if the following (1)--(4) hold for $F\in Z$.

(1) $Z=\exp(\sig_\bC)F$.

(2)  For any $N\in \sigma$, there exists $W$-relative $N$-filtration $M(N,W)$.

(3)  For any $N\in \sig$ any $p$, $NF^p\sub F^{p-1}$.

(4) If $N_{1},\dots,N_n$ generate $\sig$ and $y_j\gg0$ for any $j$, 
then $\exp(\sum_jiy_jN_j)F\in D$.
\medskip

A {\it weak fan $\Sig$ in $\fg_\bQ$} is a set of nilpotent cones in $\fg_\bR$, 
defined over $\bQ$, which satisfies the following three conditions.

(5) Every $\sig\in \Sig$ is admissible relative to $W$.

(6) If $\sig\in \Sig$ and $\tau\prec\sig$, then $\tau\in \Sig$.

(7) If $\sig,\sig'\in\Sig$ have a common interior point and 
if there exists $F\in\Dc$ such that $(\sig,F)$ and $(\sig',F)$ 
generate nilpotent orbits, then $\sig=\sig'$.
\medskip

Let $\Sig$ be a weak fan and $\G$ be a subgroup of $G_\bZ$.
$\Sig$ and $\G$ are {\it strongly compatible} 
if the following two conditions are satisfied.

(8) If $\sig\in\Sig$ and $\g\in\G$, then $\Ad(\g)\sig\in\Sig$.

(9) For any $\sig\in\Sig$, $\sig$ is generated  by $\log\G(\sig)$, 
where $\G(\sig):=\Gamma\cap\exp(\sigma)$.

\bigskip

{\bf 1.7. Moduli of log mixed Hodge structures of type $\Phi$}
\medskip

Let $\Phi:=(\Lambda, \Sig, \G)$ be a data consisting of  a Hodge data $\Lambda$ (in Notation), a weak fan $\Sig$ and a subgroup $\G$ of $G_\bZ$ such that $\Sig$ and $\G$ are strongly compatible (1.6).

Let $\sigma\in\Sigma$ and $S:=\toric_\sig$.
Let $H_\sig=(H_{\sig,\bZ},W,(\lan\;,\;\ran_w)_w)$ be the canonical local system 
$H_{\sig,\bZ}$ endowed with the weight filtration $W$ and the polarizations 
$\lan\;,\;\ran_w$ on the graded quotients $\gr^W_w$ ($w\in \bZ$) over $S^{\log}$, which is determined by the representation $\G\sub G_\bZ=\Aut(H_0,W,(\lan\;,\;\ran_w)_w)$.

Let $\Ec_\sig:=\toric_\sig\times\Dc$.
The {\it universal pre-graded polarized log mixed Hodge structure $H$ on $\Ec_\sig$} is given by $H_\sig$ together with the isomorphism 
$\cO_{\Ec_\sig}^{\log}\otimes_\bZ H_{\sig,\bZ}=\cO_{\Ec_\sig}^{\log}\otimes_{\cO_{\Ec_\sig}}H_\cO$ (1.5), where $H_\cO:=\cO_{\Ec_\sig}\otimes H_0$ is the free 
$\cO_{\Ec_\sig}$-module coming from that on $\Dc$ endowed with the universal Hodge filtration $F$.

Let 
$E_\sig:=\{x\in\Ec_\sig\,|\,\text{$H|_x$ is a graded polarized log mixed Hodge structure on $x$}\}$.
Note that slits appear in $E_\sigma$ because of log-pointwise Griffiths transversality 1.3 (2) and positivity 1.3 (3), or equivalently 1.4 (3) and 1.4 (4) respectively.

As a set, define 
$D_\Sig:=\{(\sig,Z)\,|\,\text{nilpotent orbit, $\sig\in\Sig$, $Z\sub\Dc$}\}$.
Let $\sig\in\Sig$.
Assume that $\Gamma$ is neat.
A structure as an object of $\cB(\log)$ on $\G\bs D_\Sig$ is introduced by a diagram:
$$
\CD
E_\sig@.\overset\text{GPLMH}\to\sub@.\quad\Ec:=\toric_\sig\times\Dc\\
@VV\text{$\sig_\bC$-torsor}V\\
\G(\sig)^\gp\bs D_\sig\\
@VV\text{loc.\ isom.}V\\
\G\bs D_\Sig
\endCD
$$
The action of $h\in\sig_\bC$ on $(\be(a)0_\rho,F)\in E_\sig$ is $(\be(h+a)0_\rho,\exp(-h)F)$, and the projection is $(\be(a)0_\rho,F)\mapsto[\rho,\exp(\rho_\bC+a)F]$.

Let $S\in\cB(\log)$.
A {\it log mixed Hodge structure of type $\Phi$ on $S$} is a pre-graded polarized log mixed Hodge structure $H=(H_\bZ,W,(\lan\;,\;\ran_w)_w,H_\cO)$ endowed with 
$\G$-level structure $\mu\in H^0(S^{\log},\G\bs\cI som((H_\bZ,W,(\lan\;,\;\ran_w)_w), (H_0,W,(\lan\;,\;\ran_w)_w)))$ satisfying the following condition:
For any point $s\in S$, any point $t\in\tau^{-1}(s)=s^{\log}$ and any representative $\tilde\mu_t:H_{\bZ,t}\tra H_0$, there exists $\sig\in\Sig$ such that $\sig$ contains $\tilde\mu_t\pi_1^+(s^{\log})\tilde\mu_t^{-1}$ and $(\sig,\tilde\mu_t(\bC\otimes_{\cO_{S,t}^{\log}}F_t))$ 
generates a nilpotent orbit.
Here $\pi_1^+(s^{\log}):=\operatorname{Image}(\Hom((M_S/\cO_S^\times)_s,\bN)\hookrightarrow\pi_1(s^{\log})\to\Aut(H_{\bZ,t}))$ is the local monodromy monoid of $H_\bZ$ at $s$ (cf.\ [KU09, 3.3.2]).
(Then, the smallest such $\sig$ exists.)

\proclaim{Theorem}
For a given data $\Phi$, we have the following.

\noindent
{\rm(i)} $\G\bs D_\Sig\in\cB(\log)$, which is Hausdorff.
If $\G$ is neat, $\G\bs D_\Sig$ is a log manifold.

\noindent
{\rm(ii)} On $\cB(\log)$, $\G\bs D_\Sig$ represents a functor 
$\LMH_\Phi$ of log mixed Hodge structures of type $\Phi$.
\endproclaim

\proclaim{Log period map}
Given $\Phi$.
Let $S\in\cB(\log)$.
Then we have an isomorphism
$$\LMH_\Phi(S)\tra\Map(S,\G\bs D_\Sig), \;
H\mapsto\big(S\ni s\mapsto[\sig,\exp(\sig_\bC)\tilde\mu_t(\bC\otimes_{\cO_{S,t}^{\log}}F_t)])\;\;
(t\in\tau^{-1}(s)),
$$
which is functorial in $S$.
\endproclaim

A log period map is a unified compactification of a period map and a normal function 
of Griffiths.

The above $\G\bs D_\Sig$ is the fine moduli of log mixed Hodge structures of type $\Phi$, whose underlying coarse moduli, in the sense of log points, is the set of equivalence classes of all nilpotent orbits of specified type.

\bigskip

{\bf 1.8. N\'eron model for admissible normal function}
\medskip

We review some results from [KNU14, Theorem 1.3], [KNU13, III, Section 6.1] and [KNU10, Section 8] adapted to the situation (5) in Introduction.

For a pure case $h^{p,q}=1$ ($p+q=3$, $p,q\ge0$) and $h^{p,q}=0$ otherwise, a complete fan is constructed in [KU09, Section 12.3].
For a mixed case $h^{p,q}=1$ (the above $(p,q)$, plus $(p,q)=(2,2)$) and $h^{p,q}=0$ otherwise, over the above fan, a fan of N\'eron model for given admissible normal function is constructed in [KNU14, Theorem 3.1], and we have a N\'eron model in the following sense.

Let $S\in\cB(\log)$, $U:=S_{\triv}\subset S$ (consisting of those points with trivial log structure), $H_{(-1)}$ be a polarized variation of Hodge structure of weight $-1$ (Tate-twisted by 2 for $\cH$ in Introduction (5)) on $U$ and $L_\bQ$ be a local system of $\bQ$-vector spaces which is an extension of $\bQ$ by $H_{(-1),\bQ}$.
An admissible normal function over $U$ for $H_{(-1)}$ underlain by the local system $L_\bQ$ can be regarded as an admissible variation of mixed Hodge structure which is an extension of $\bZ$ by $H_{(-1)}$ and lies over local system $L_\bQ$.

For any given unipotent admissible normal function over $U$ as above, $H_{(-1)}$ and $L_\bQ$ extend to a polarized log mixed Hodge structure on $S$ and a local system on $S^{\log}$, respectively, denoted by the same symbols, and there is a relative log manifold $J_{L_{\bQ}}$ over $S$ which is strict over $S$ (i.e., endowed with the pullback log structure from $S$) and which represents the following functor on $\cB/S^{\circ}$ ($S^{\circ}\in\cB$ is the underlying space of $S$):

$S'\mapsto$ \{LMH $H$ on $S'$ satisfying $H(\gr^W_w)=H_{(w)}|_{S'}$ $(w=-1,0)$ and ($\ast$) below\}/isom.

($\ast$) Locally on $S'$, there is an isomorphism $H_{\bQ} \simeq L_{\bQ}$  on $(S')^{\log}$ preserving  $W$.
\medskip

\noindent
Here $H_{(w)}|_{S'}$ is the pullback of $H_{(w)}$ by the structure morphism $S' \to S^{\circ}$, 
and $S'$ is endowed with the pullback log structure from $S$.

Put $H':=H_{(-1)}$.
In the present case, we have $J_{L_\bQ}=\cExt^1_{\text{LMH}/S}(\bZ, H')$ by [KNU13, Corollary 6.1.6].
This is the subgroup of $\tau_*(H'_{\cO^{\log}}/(F^0+H'_\bZ))$ restricted by admissibility condition and log-pointwise Griffiths transversality condition ([KNU10, Section 8], cf.\ 1.5).
Let $\tilde J_{L_\bQ}$ be the pullback of $J_{L_\bQ}$ by $\tau_*(H'_{\cO^{\log}}/F^0)\to \tau_*(H'_{\cO^{\log}}/(F^0+H'_\bZ))$, and $\bar J_{L_\bQ}$ be the image of $\tilde J_{L_\bQ}$ by $\tau_*(H'_{\cO^{\log}}/F^0)\to \tau_*(H'_{\cO^{\log}}/F^{-1})$.
Then, by using the polarization, we have a commutative diagram:
$$
\matrix
J_{L_\bQ}&=&\cExt^1_{\text{LMH}/S}(\bZ, H')&\subset&
\tau_*(H'_{\cO^{\log}}/(F^0+H'_\bZ))&@>\text{pol}>\sim>&
\tau_*((F^0)^*/H'_\bZ) \\
@AAA@.@AAA@AAA\\
\tilde J_{L_\bQ}&&\subset&&
\tau_*(H'_{\cO^{\log}}/F^0)&@>\text{pol}>\sim>&
\tau_*((F^0)^*) \\
@VVV@.@VVV@VVV\\
\bar J_{L_\bQ}&&\subset&&
\tau_*(H'_{\cO^{\log}}/F^{-1})&@>\text{pol}>\sim>&\;\tau_*((F^1)^*).
\endmatrix
$$

\bigskip

{\bf 2. Quintic threefolds}
\medskip

Let $V$ be a quintic threefold in $\bP^4$ and let $V_\psi^\circ$ be its mirror family (cf.\ [CoK99, Sect.\ 4.2]).

In this section, we give a correspondence table of A-model for $V$ and B-model for $V_\psi^\circ$.
This is a precision and an expansion of our previous [U14, 3] by using $\hat\Gamma$-integral structure of Iritani [I11].
We will use this table in Section 3 below.
\bigskip

{\bf 2.1. Quintic mirror family}
\medskip

Following [M93], [MW09], etc., we briefly recall the construction of the mirror family $V_\psi^\circ$ by quotient method.
Let $V_\psi:  f:=\sum_{j=1}^5x_j^5-5\psi\prod_{j=1}^5 x_j=0$ 
$(\psi\in\bP^1)$ be the Dwork pencil of quintics in $\bP^4$.
Let $\mu_5$ be the group consisting of the fifth roots of the unity in $\bC$.
Then the group $G:=\{(a_j)\in(\mu_5)^5\,|\,a_1\dots a_5=1\}$ acts on $V_\psi$ by $x_j\mapsto a_jx_j$.
Let $V_\psi^\circ$ be a crepant resolution of quotient singularity of $V_\psi/G$ (cf.\ [M93], [MW09]).
Divide further by the action $(x_1,\dots, x_5)\mapsto (a^{-1}x_1,x_2,\dots, x_5)$ and $\psi\mapsto a\psi$ $(a\in\mu_5)$.

\bigskip

{\bf 2.2. Picard--Fuchs equation on the mirror $V^\circ$}
\medskip

Let $\Omega$ be a 3-form on $V_\psi^\circ$ with a log pole over $\psi=\infty$ induced from 
$$
\Big(\frac{5}{2\pi i}\Big)^3\Res_{V_\psi}\Big(\frac{\psi}{f}\sum_{j=1}^5(-1)^{j-1}x_jdx_1\wedge\dots\wedge\widehat{dx_j}\wedge\dots\wedge dx_5\Big).
$$
Let $z:=1/(5\psi)^5$ and $\theta:=zd/dz$.
Let 
$$
\cL:=\theta^4-5z(5\theta+1)(5\theta+2)(5\theta+3)(5\theta+4)
$$
be the Picard--Fuchs differential operator for $\Omega$, i.e., $\cL\Omega=0$ via 
the Gauss-Manin connection $\nabla$.
There are three special points of the complex moduli:

$z=0$ : maximally unipotent monodromy point, 

$z=\infty$ : Gepner point, 

$z=1/5^5$ : conifold point.

At $z=0$, the Picard--Fuchs differential equation $\cL y=0$ has the indicial equation 
$\rho^4=0$ ($\rho$ is indeterminate), i.e., maximally unipotent.
By the Frobenius method, we have a basis of solutions $y_j(z)$ $(0\le j\le3)$ as follows.
Let
$$
\t y(-z;\rho):=\sum_{n=0}^\infty\frac{\prod_{m=1}^{5n}(5\rho+m)}{\prod_{m=1}^n(\rho+m)^5}(-z)^{n+\rho}
$$
be a solution of $\cL(\t y(-z;\rho))=\rho^4(-z)^\rho$, and let
$$
\t y(-z;\rho)=y_0(z)+y_1(z)\rho+y_2(z)\rho^2+y_3(z)\rho^3+\cdots,\quad
y_j(z):=\frac{1}{j!}\frac{\partial^i\t y(-z;\rho)}{\partial\rho^j}|_{\rho=0}
$$
be the Taylor expansion at $\rho=0$.
Then, $y_j$ $(0\le j\le 3)$ form a basis of homogeneous solutions for the linear differential equation $\cL y=0$.
We have 
$$
\align
&y_0=f_0=\sum_{n=0}^\infty\frac{(5n)!}{(n!)^5}z^n,\\
&y_1=f_0\log z+f_1=y_0\log z+5\sum_{n=1}^\infty\frac{(5n)!}{(n!)^5}\Big(\sum_{j=n+1}^{5n}\frac{1}{j}\Big)z^n,\\
&2!y_2=f_0(\log z)^2+2f_1\log z+f_2,\\
&3!y_3=f_0(\log z)^3+3f_1(\log z)^2+3f_2\log z+f_3,
\endalign
$$
where all $f_j$ are holomorphic functions in $z$ with $f_0(0)=1$ and $f_j(0)=0$ for $j>0$.

Define the canonical parameters by $t:=y_1/y_0$, $u:=t/2\pi i$, and the canonical coordinate by $q:=e^t=e^{2\pi iu}$ which is a specific chart of the log structure given by the divisor $(z=0)$ of a disc in $\bP^1$ and gives a mirror map.

Write $z=z(q)$ which is holomorphic in $q$.
Then we have
$$
\log z=2\pi iu-\frac{5}{y_0(z(q))}\sum_{n=1}^\infty\frac{(5n)!}{(n!)^5}\Big(\sum_{j=n+1}^{5n}\frac{1}{j}\Big)z(q)^n.
$$

The Gauss-Manin potential of $V_z^\circ$ is 
$$
\Phi_{\text{\rm GM}}^{V^\circ}
=\frac{5}{2}\Big(\frac{y_1}{y_0}\frac{y_2}{y_0}-\frac{y_3}{y_0}\Big).
$$

Let $\tilde\Omega:=\Omega/y_0$ and $\delta:=2\pi iqd/dq=du$.
Then, the Yukawa coupling at $z=0$ is 
$$
Y:=-\int_{V^\circ}\tilde\Omega\wedge\nabla_\delta\nabla_\delta\nabla_\delta\tilde\Omega=\frac{5}{(1+5^5z)y_0(z)^2}\Big(\frac{qdz}{zdq}\Big)^3.
$$

\bigskip

{\bf 2.3. A-model of quintic $V$}
\medskip

Let $V$ be a general quintic hypersurface in $\bP^4$.
Let $H$ be the cohomology class of a hyperplane section of $V$ in $\bP^4$, $K(V)=\bR_{>0}H$ be the K\"ahler cone of $V$, and $u$ be the coordinate of $\bC H$. 
Put $t:=2\pi iu$.
A complexified K\"ahler moduli is defined as 
$$
\cK\cM(V):=(H^2(V,\bR)+iK(V))/H^2(V,\bZ) \tra \Delta^*,\quad
uH\mapsto q:=e^{2\pi iu}.
$$
Let $C\in H_2(V,\bZ)$ be the homology class of a line on $V$.

For $\b=dC\in H_2(V,\bZ)$, define $q^\b:=q^d$.
The Gromov--Witten potential of $V$ is defined as 
$$
\Phi_{\text{\rm GW}}^V:=\frac{1}{6}\int_V(2\pi iuH)^3
+ \sum_{0\ne\b\in H_2(V,\bZ)}N_dq^\b
=\frac{5}{6}(2\pi i)^3u^3 + \sum_{d>0}N_dq^d.
$$
Here the Gromov--Witten invariant $N_d$ is
$$
\align
&\overline M_{0,0}(\bP^4,d)@<\pi_1<<
\overline M_{0,1}(\bP^4,d)@>e_1>>\bP^4,\\
&N_d:=\int_{\overline M_{0,0}(\bP^4,d)}c_{5d+1}(\pi_{1*}e_1^*\cO_{\bP^4}(5)).
\endalign
$$
Note that $N_d=0$ if $d\le0$.
Let $N_d=\sum_{k|d}n_{d/k}k^{-3}$.
Then $n_{d/k}$ is the instanton number.
($n_{l}$ here is different from $n_{l}$ in (4) in Introduction.)

The differentials of $\Phi=\Phi_{\text{\rm GW}}^V$ are computed easily:
$$
\frac{d\Phi}{du}=\frac{5}{2}(2\pi i)^3u^2 + (2\pi i)\sum_{d>0}N_ddq^d,\quad
\frac{d^2\Phi}{du^2}=5(2\pi i)^3u + (2\pi i)^2\sum_{d>0}N_dd^2q^d.
$$

\bigskip

{\bf 2.4. Integral structure}
\medskip

As we stated in Introduction, we consider everything near the large radius point $q_{0}$ and the maximally unipotent monodromy point $p_{0}$.
Let $S$ be a neighborhood disc of $q_0$ (resp.\ $p_0$) in $\overline{\Cal K\Cal M}(V)$ (resp.\ $\overline\cM(V^\circ)$) for A-model of $V$ (resp.\ for B-model of $V^\circ$), and let $S^*$ be $S\smallsetminus \{q_0\}$ (resp.\ $S\smallsetminus \{p_0\}$) for A-model (resp.\ B-model) (see 2.2, 2.3).
Endow $S$ with the log structure associated to the divisor $S\smallsetminus S^*$.

The B-model variation of Hodge structure $\cH^{V^\circ}$ is the usual variation of Hodge structure arising from the smooth projective family $f:X\to S^*$ of the quintic mirrors over the punctured neighborhood of $p_0$.
Its integral structure is the usual one $\cH^{V^\circ}_\bZ=R^3f_*\bZ$.
This is compatible with the monodromy weight filtration $M$ around $p_0$.
Define $M_{k,\bZ}:=M_k\cap \cH^{V^\circ}_\bZ$ for all $k$.

For the A-model $\cH^V$ on $S^*$, the locally free sheaf on $S^*$, the Hodge filtration, and the monodromy weight filtration $M$ around $q_0$ are given by $\cH^V_\cO:=\cO_{S^*}\otimes (\bigoplus_{0\le p\le 3}H^{2p}(V))$, $F^p:=\cO_{S^*}\otimes H^{\le2(3-p)}(V)$, and $M_{2p}:=H^{\ge2(3-p)}(V)$, respectively.
Iritani defined $\hat\Gamma$-integral structure in more general setting in [I11, Definition 3.6].
In the present case, it is characterized as follows.
Let $H$ and $C$ be a hyperplane section and a line on $V$, respectively.
Then, in the present case, a basis of the $\hat\Gamma$-integral structure is given by $\{s(\cE)\,|\, \cE\;\text{is}\; \cO_V,\cO_H,\cO_C,\cO_{\text{pt}}\}$ [ibid, Example 6.18], where $s(\cE)$ is a unique $\nabla^{\text{even}}$-flat section satisfying an asymptotic condition
$$
s(\cE)\sim(2\pi i)^{-3}e^{-2\pi iuH}\cdot \hat\Gamma(T_V)\cdot(2\pi i)^{\deg/2}\operatorname{ch}(\cE)
$$
at the large radius point $q_0$ when $\Im(u)\to\infty$ for each fixed $\Re(u)$.
Here, for the Chern roots $c(T_V)=\prod_{j=1}^3(1+\delta_j)$, the Gamma class $\hat\Gamma(T_V)$ is defined by 
$$
\align
\hat\Gamma(T_V):=\prod_{j=1}^3\Gamma(1+\delta_j)&=\exp(-\gamma c_1(V)+\sum_{k\ge2}(-1)^k(k-1)!\zeta(k)\operatorname{ch}_k(T_V))
\\
&=\exp(\zeta(2)\operatorname{ch}_2(T_V)-2\zeta(3)\operatorname{ch}_3(T_V)) 
\endalign
$$
where $\gamma$ is the Euler constant, and $\deg|_{H^{2p}(V)}:=2p$.
The important point is that this class $\hat\Gamma(T_V)$ plays the role of a \lq\lq square root" of the Todd class in Hirzebruch-Riemann-Roch  ([I09, 1], [I11, 1, (13)]).
Denote this $\hat\Gamma$-integral structure by $\cH^V_\bZ$.
This is compatible with the monodromy weight filtration $M$ and we define $M_{k,\bZ}:=M_k\cap \cH^V_\bZ$ for all $k$.
The above asymptotic relation is actually computed as
$$
\align
s^0&:=s(\cO_{\text{pt}})=\frac{1}{5}H^3,\\
s^1&:=s(\cO_C)=\frac{1}{5}(2\pi i)^{-1}H^2+\frac{1}{5}(-u+1)H^3,\\
s^2&:=s(\cO_H)\sim (2\pi i)^{-2}H+\frac{5}{2}(2\pi i)^{-1}\Big(-u-\frac{1}{2}\Big)H^2+\Big(\frac{1}{2}u^2+\frac{1}{2}u+\frac{7}{12}\Big)H^3,\\
s^3&:=s(\cO_V)\\
&\sim (2\pi i)^{-3}-(2\pi i)^{-2}uH+(2\pi i)^{-1}\Big(\frac{1}{2}u^2+\frac{5}{12}\Big)H^2+\Big(-\frac{1}{6}u^3-\frac{5}{12}u+\frac{5i\zeta(3)}{\pi^3}\Big)H^3.
\endalign
$$

Fixing an isomorphism of VHS in (2) in Introduction, we also use $s^p$ for the corresponding $\nabla$-flat integral basis for the B-model $\cH^{V^\circ}_\bZ$ (vanishing cycles are used for B-model in [I11, Theorems 6.9, 6.10, Example 6.18]).

In both A-model case and B-model case, the integral structures $\cH^V_\bZ$ and 
$\cH^{V^\circ}_\bZ$ on $S^*$ extend to the local systems of $\bZ$-modules over $S^{\log}$ ([O03], [KU09, Proposition 2.3.5]), still denoted $\cH^V_\bZ$ and 
$\cH^{V^\circ}_\bZ$, respectively.

Consider a diagram:
$$
\CD
\tilde S^{\log}:=\bR\times i(0,\infty]@.\;\supset@.\;\tilde S^*:=\bR\times i(0,\infty)\\
@VVV@.@VVV\\
S^{\log}@.\;\supset @.S^*\\
@V{\tau}VV\\
S
\endCD
$$
The coordinate $u$ of $\tilde S^*$ extends over $\tilde S^{\log}$.
Fix base points as $u_0=0+i\infty\in\tilde S^{\log}\mapsto b:=\bar 0+i\infty\in S^{\log}\mapsto q=0\in S$,
where $q=0$ corresponds to $q_0$ for A-model and $p_0$ for B-model.
Note that fixing a base point $u=u_0$ on $\tilde S^{\log}$ is equivalent  to fixing
a base point $b$ on $S^{\log}$ and also a branch of $(2\pi i)^{-1}\log q$.

Let  $B:=\cH^V_\bZ(u_0)=\cH^V_\bZ(b)$ for A-model and $B:=\cH^{V^\circ}_\bZ(u_0)=\cH^{V^\circ}_\bZ(b)$ for B-model.

\bigskip


{\bf 2.5. Correspondence table}
\medskip

We use the mirror theorems (1)--(3) in Introduction.
Put $\Phi:=\Phi_{\text{GW}}^V=\Phi_{\text{GM}}^{V^\circ}$ and fix an isomorphisim of VHS in (2) in Introduction (cf.\ 2.4).
\medskip

(1A) {\it Polarization of A-model of $V$.}
$$
S(\a,\b):=(-1)^p(2\pi i)^3\int_V\a\cup\b\quad
(\a\in H^{p,p}(V), \b\in H^{3-p,3-p}(V)).
$$
\medskip

(1B) {\it Polarization of B-model of $V^\circ$.}
$$
Q(\a,\b):=(-1)^{3(3-1)/2}\int_{V^\circ}\a\cup\b
=-\int_{V^\circ}\a\cup\b\quad
(\a,\b\in H^3(V^\circ)).
$$
\medskip

(2A) {\it $\bZ$-basis compatible with monodromy weight filtration.}

Let  $B:=\cH^V_\bZ(u_0)=\cH^V_\bZ(b)$ be as in Section 2.4.
Let $b^3:=s^3(u_0)=s(\cO_V)(u_0)$, $b^2:=s^2(u_0)=s(\cO_H)(u_0)$, $b^1:=s^1(u_0)=s(\cO_C)(u_0)$ and $b^0:=s^0(u_0)=s(\cO_{\text{pt}})(u_0)$ be the basis of the fiber $B$ at $u_0$ coming from $\nabla$-flat integral basis in  2.4.

The endomorphism of $B_\bQ:=\bQ\otimes B$ coming from the monodromy logarithm coincides with the cup product with $-2\pi iH$ where $H$ is a hyperplane section of $V$ ([I11, Definition 3.6], cf.\ the third Proposition in 1.4).
Hence the above basis is compatible with the monodromy weight filtration $M$.

\medskip

(2B) {\it $\bZ$-basis compatible with monodromy weight filtration.}

Let $B:=\cH^{V^\circ}_\bZ(u_0)=\cH^{V^\circ}_\bZ(b)$, and $b^0, b^1, b^2, b^3$ be the basis of $B$ corresponding to that in (2A) by the mirror symmetry (2) and (3) in Introduction.

The endomorphism of $B_\bQ$ coming from the monodromy logarithm is denoted by $N$, and the above basis is compatible with the monodromy weight filtration $M$ [ibid].

\medskip

For both cases (2A) and (2B), we regard $B$ as a constant sheaf on $S^{\log}$ and also on $S$, endowed with the associated filtrations $M$.

\medskip

From the asymptotics of the basis $s^p$ $(0\le p\le3)$ in 2.4, the matrix of the polarization pairings $S$ in (1A) and $Q$ in (1B) for the basis  $b^p=s^p(u_0)$ is computed as
$$
(S(b^p,b^q))_{p,q}=(Q(b^p,b^q))_{p,q}=
\pmatrix
0&0&0&-1\\
0&0&1&-1\\
0&-1&0&-5\\
1&1&5&0
\endpmatrix.
$$

\medskip

(3A) {\it Sections compatible with Deligne decomposition and inducing $\bZ$-basis of $\gr^M$ for A-model of $V$.}

Let $T^3$, $T^2$, $T^1$, and $T^0$ be the basis of $\cH^V_{\cO}$ corresponding to the  $e^3$, $e^2$, $e^1$, and $e^0$ in (3B) below by the mirror symmetry (2) and (3) in Introduction. 
Then $S(T^3,T^0)=1$ and $S(T^2,T^1)=-1$.
Hence $T^3$, $T^2$, $-T^0$, $T^1$ form a symplectic base for $S$ in (1A).

Note that on $\gr^M$ they are 
$$
\align
&\gr^M_3(T^3)=1\in H^0(V,\bZ),\quad
\gr^M_2(T^2)=H\in H^2(V,\bZ),\\
&\gr^M_1(T^1)=C\in H^4(V,\bZ),\quad
\gr^M_0(T^0)=[\text{pt}]\in H^6(V,\bZ),
\endalign
$$
where $H$ and $C$ are the cohomology classes of a hyperplane section and a line on $V$, respectively.
Abusing notation, we mean by $C$ the Poincar\'e dual class of the homology class in 2.3.

\medskip

(3B) {\it Sections compatible with Deligne decomposition and inducing $\bZ$-basis of $\gr^M$ for B-model of $V^\circ$.}

We use Deligne decomposition [D97].
We consider $B$ in (2B) as a constant sheaf on $S^{\log}$.
We have locally free $\cO_S$-submodules $\cM_{2p}:=\tau_*(\cO_S^{\log}\otimes_\bZ M_{2p}B)$ and $\cF^p$ in $\tau_*(\cO_S^{\log}\otimes_\bZ B)=\cO_S\otimes_\bZ B$ (canonical extension of Deligne in the second Proposition in 1.4).
The mixed Hodge structure of Hodge--Tate type $(\cM,\cF)$ has decomposition:
$$
\cO_S\otimes_\bZ B=\bigoplus_pI^{p,p}, \qquad
I^{p,p}:=\cM_{2p}\cap\cF^p\tra \gr^{\cM}_{2p}.
$$
Transporting the basis $b^p$ $(0\le p\le3)$ of $B$ in (2B), regarded as sections of the constant sheaf $B$ on $S^{\log}$, via isomorphism
$$
I^{p,p}@>\sim>>\cO_S\otimes_\bZ \gr^M_{2p}B
$$
we define sections $e^p\in I^{p,p}$ $(0\le p\le3)$ over $S$.
Then, $e^3$, $e^2$,$-e^0$, $e^1$ form a symplectic basis for $Q$ in (1B), and 
$e^3=\tilde\Omega=\Omega/y_0$ over $S$.
\medskip

The asymptotic relation of the $T^p$ in (3A) (resp.\ the $e^p$ in (3B)) can be computed, via the $s^p$, from (7A) (resp.\ (7B)) below.

\medskip

(4A) {\it A-model connection $\nabla=\nabla^{\text{even}}$ of $V$.}

Let $\delta=d/du=2\pi iqd/dq$.
The Dubrovin connection $\nabla$ (cf.\ [CK91, Sect.\ 8.4]) is characterized by 
$$
\align
&\nabla_\delta T^0=0,\quad
\nabla_\delta T^1=T^0,\\
&\nabla_\delta T^2=\frac{1}{(2\pi i)^3}\frac{d^3\Phi}{du^3}T^1
=\Big(5+\frac{1}{(2\pi i)^3}\frac{d^3\Phi_{\text{hol}}}{du^3}\Big)T^1,\\
&\nabla_\delta T^3=T^2.
\endalign
$$
$\nabla$ is flat, i.e., $\nabla^2=0$, and extends to a log connection over $S^{\log}$.
\medskip

(4B) {\it B-model connection $\nabla=\nabla^{\operatorname{GM}}$ of $V^\circ$.}

Let $\delta=d/du=2\pi iqd/dq$.
The Gauss-Manin connection $\nabla$ is computed as 
$$
\align
&\nabla_\delta e^0=0,\quad
\nabla_\delta e^1=e^0,\\
&\nabla_\delta e^2
=\frac{1}{(2\pi i)^3}\frac{d^3\Phi}{du^3}e^1
=Ye^1
=\frac{5}{(1+5^5)y_0(z)^2}
\Big(\frac{q}{z}\frac{dz}{dq}\Big)^3e^1,\\
&\nabla_\delta e^3=e^2.
\endalign
$$
$\nabla$ is flat, i.e., $\nabla^2=0$, and extends to a log connection over $S^{\log}$.

\medskip


(5A) {\it $\nabla$-flat basis of  $\cH_\bC^V$ inducing $\bZ$-basis of $\gr^M$.}
$$
\align
\tilde s^0&:=T^0,\\
\tilde s^1&:=T^1-uT^0,\\
\tilde s^2&:=T^2-\frac{1}{(2\pi i)^3}\frac{d^2\Phi}{d u^2}T^1+\frac{1}{(2\pi i)^3}\frac{d\Phi}{d u}T^0,\\
\tilde s^3&:=T^3-uT^2
+\frac{1}{(2\pi i)^3}
\Big(u\frac{d^2\Phi}{d u^2}-\frac{d\Phi}{dl u}\Big)T^1
-\frac{1}{(2\pi i)^3}\Big(u\frac{d\Phi}{d u}-2\Phi\Big)T^0.
\endalign
$$
Then $\tilde s^3$, $\tilde s^2$,$-\tilde s^0$, $\tilde s^1$ form a symplectic basis for $S$ in (1A).

\medskip

(5B) {\it $\nabla$-flat basis for $\cH_\bC^{V^\circ}$ inducing $\bZ$-basis of $\gr^M$.}
$$
\align
\tilde s^0&:=e^0,\\
\tilde s^1&:=e^1-ue^0,\\
\tilde s^2&:=e^2-\frac{1}{(2\pi i)^3}\frac{d^2\Phi}{d u^2}e^1+\frac{1}{(2\pi i)^3}\frac{d \Phi}{d u}e^0,\\
\tilde s^3&:=e^3-ue^2
+\frac{1}{(2\pi i)^3}
\Big(u\frac{d^2\Phi}{d u^2}-\frac{d\Phi}{d u}\Big)e^1
-\frac{1}{(2\pi i)^3}\Big(u\frac{d\Phi}{d u}-2\Phi\Big)e^0.
\endalign
$$
Then $\tilde s^3$, $\tilde s^2$,$-\tilde s^0$, $\tilde s^1$ form a symplectic basis for $Q$ in (1B).

\medskip

For both cases (5A) and (5B), by using (4A) and (4B), the $\nabla$-flat bases  $\tilde s^p$ are determined inductively on $0\le p\le3$ from the $T^p$ in (3A) and the $e^p$ in (3B).
These $\nabla$-flat bases $\tilde s^p$ are characterized by the Frobenius solutions $y_j$ $(0\le j \le 3)$ in 2.2 such that $y_0T^3$ in A-model and $y_0e^3=\Omega$ in B-model coincide with   
$$
y_0\tilde s^3+(2\pi i)^{-1}y_1\tilde s^2+5(2\pi i)^{-2}y_2\tilde s^1+5(2\pi i)^{-3}y_3\tilde s^0.
$$


\medskip


(6A), (6B) {\it Relations of $\nabla$-flat $\bZ$-basis $s^p$ and the $\nabla$-flat basis $\tilde s^p$.}
$$
\align
s^0&=\tilde s^0,\\
s^1&=\tilde s^1+\tilde s^0,\\
s^2&=\tilde s^2-\frac{5}{2}\tilde s^1+\frac{35}{12}\tilde s^0,\\
s^3&=\tilde s^3+\frac{25}{12}\tilde s^1+\frac{25i\zeta(3)}{\pi^3}\tilde s^0.
\endalign
$$

\medskip


(7A) {\it Expression of the $T^p$ by the $s^p$ over $S^{\log}$.}

It is computed that $T^p$ are written by the $\nabla$-flat $\bZ$-basis $s^p$ of $\cH_\bZ^V$ as follows.
$$
\align
T^0&=s^0,\\
T^1&=s^1+(u-1)s^0,\\
T^2&=s^2+\Big(\frac{1}{(2\pi i)^3}\frac{d^2\Phi}{d u^2}+\frac{5}{2}\Big)s^1
+\Big(\frac{1}{(2\pi i)^3}\Big(u\frac{d^2\Phi}{d u^2}-
\frac{d \Phi}{d u}\Big)-5u-\frac{65}{12}\Big)s^0,\\
T^3&=s^3+us^2
+\Big(\frac{1}{(2\pi i)^3}\frac{d\Phi}{d u}
+\frac{5}{2}u-\frac{25}{12}\Big)s^1\\
&\quad
+\Big(
\frac{1}{(2\pi i)^3}
\Big(u\frac{d\Phi}{d u}-2\Phi\Big)
-\frac{65}{12}u+\frac{25}{12}-\frac{25i}{\pi^3}\zeta(3)\Big)s^0.
\endalign
$$
\medskip

(7B) {\it Expression of the $e^p$ by the $s^p$ over $S^{\log}$.}

It is computed that $e^p$ are written by the $\nabla$-flat $\bZ$-basis $s^p$ of $\cH_\bZ^{V^\circ}$ as follows.
$$
\align
e^0&=s^0,\\
e^1&=s^1+(u-1)s^0,\\
e^2&=s^2+\Big(\frac{1}{(2\pi i)^3}\frac{d^2\Phi}{d u^2}+\frac{5}{2}\Big)s^1
+\Big(\frac{1}{(2\pi i)^3}\Big(u\frac{d^2\Phi}{d u^2}-
\frac{d \Phi}{d u}\Big)-5u-\frac{65}{12}\Big)s^0,\\
e^3&=s^3+us^2
+\Big(\frac{1}{(2\pi i)^3}\frac{d\Phi}{d u}
+\frac{5}{2}u-\frac{25}{12}\Big)s^1\\
&\quad
+\Big(
\frac{1}{(2\pi i)^3}
\Big(u\frac{d\Phi}{d u}-2\Phi\Big)
-\frac{65}{12}u+\frac{25}{12}-\frac{25i}{\pi^3}\zeta(3)\Big)s^0.
\endalign
$$

\medskip

(8A), (8B) {\it Relations of integral periods and Frobenius solutions.}

Let $\eta_j$ $(0\le j \le 3)$ be the integral periods defined by the condition that $y_0T^3$ in A-model and $y_0e^3=\Omega$ in B-model coincide with 
$\eta_0s^3+\eta_1s^2+\eta_2s^1+\eta_3s^0$.
Then the relations in (6A), (6B) are interpreted as 
$$
\align
\eta_0&=y_0,\\
\eta_1&=(2\pi i)^{-1}y_1,\\
\eta_2&=5(2\pi i)^{-2}y_2+\frac{5}{2}(2\pi i)^{-1}y_1-\frac{25}{12}y_0,\\
\eta_3&=5(2\pi i)^{-3}y_3-5(2\pi i)^{-2}y_2
+\frac{65}{12}(2\pi i)^{-1}y_1+\Big(\frac{25}{12}-\frac{25i\zeta(3)}{\pi^3}\Big)y_0.
\endalign
$$
\medskip

{\it Remark.}
The $\eta_j$ coincide with the corresponding coefficients of the expression of $y_0T^3$ in (7A) and of $y_0e^3=\Omega$ in (7B), and yield the same integral structure for periods given in [CDGP91] and [LW12].
Indeed, for the notation $\varpi_j$ in [LW12, (8)], they are related as
$\varpi_0=\eta_0$, $\varpi_1=\eta_1$, $\varpi_2=\eta_2-5\eta_1$, and $\varpi_3=-\eta_3-\eta_2-5\eta_1$.

\bigskip


{\bf 2.6. Proofs of results in 2.5}
\medskip

We prove the correspondence table in 2.5.
\medskip

(1A), (1B), (2A), (2B), (3A), and (3B) are definitions, constructions or almost direct consequences.

\medskip

{\it Proofs of (4A) and (4B) in 2.5.}
We prove (4B).
(4A) follows by  mirror symmetry theorems (1)--(3) in Introduction.

We improve the proof of [CoK99, Prop.\ 5.6.1] carefully by a log Hodge theoretic understanding in 1.4 of the relation among a constant sheaf and the local system on $S^{\log}$,  the canonical extension of Deligne on $S$, and the Deligne decomposition.

We investigate the Gauss-Manin connection $\nabla$, corresponding to the local system $\cH^{V^\circ}_\bZ$, contracted with $\delta=2\pi i qd/dq$.
Since $e^p$ maps to a $\nabla(\gr^{\cM}_{2p})$-flat element of $\gr^{\cM}_{2p}$, $\nabla_\delta(e^p)$ lies in $\cM_{2p-1}=\cM_{2p-2}$.
But $e^p$ is also an element of $\cF^p$, so that $\nabla_\delta(e^p)$ lies in $\cF^{p-1}$ by Griffiths transversality.
This shows that $\nabla_\delta(e^p)$ is an element of $I^{p-1,p-1}$, and it follows that
$$
\nabla_\delta(e^3)=Y_3e^2,\quad 
\nabla_\delta(e^2)=Y_2e^1,\quad 
\nabla_\delta(e^1)=Y_1e^0,\quad 
\nabla_\delta(e^0)=0
$$
for some $Y_1,Y_2,Y_3\in\cO^{\log}_S$.
However, since $Q(e^3,e^1)=0$ by orthogonality of Hodge filtration, we have
$$
0=\delta Q(e^3,e^1)=Q(\nabla_\delta(e^3),e^1)+Q(e^3,\nabla_\delta(e^1))
=Y_3Q(e^2,e^1)+Y_1Q(e^3,e^0)=-Y_3+Y_1,
$$
where the last equality follows from 2.5 (3B).

Since $\nabla$ has a regular singular point and $\delta=2\pi i qd/dq$, $Y_1$ is holomorphic over $S$.
Considering over the log point $p_0$, we claim $Y_1(0)=\pm1$.
Since $e^1$ is taken to be the canonical extension over $p_0$, we have $\nabla_\delta(e^1)=N(e^1)=-e^0$ by [KU09, Prop.\ 2.3.4 (ii)] (cf.\ the third proposition in 1.4).
Replacing $e^1, e^2$ by $-e^1, -e^2$, we have $Y_{1}(0)=1$.

Since we use the canonical coordinate $q$ in 2.2, the arguments in [CoK99, Sect.\ 5.6.4, Sect.\ 2.3] yield 
$$
q=\exp\Big(\int Y_1(q)\frac{dq}{q}\Big).
$$
Taking logarithm of both sides and differentiating them by $d/d\log q$, we have $Y_1(q)=1$, hence $\nabla_\delta e^1=e^0$ and $\nabla_\delta e^3=e^2$.
Thus, relative to the basis $e^0,e^1,e^2,e^3$ and using the canonical coordinate $q$, $\nabla_\delta$ has the connection matrix
$$
\pmatrix
0&1&&\\
&0&Y&\\
&&0&1\\
&&&0
\endpmatrix
$$
where $Y=Y_2$.

The proofs in [ibid, Proof of Prop.\ 5.6.1] for the following assertions work well:
$e^3$ is the normalized 3-form $\tilde \Omega$; the Picard--Fuchs equation for $\tilde \Omega$ is $\nabla_\delta^2(\nabla_\delta^2\tilde\Omega/Y)=0$; $Y$ is the Yukawa coupling.
The notation $e^0,e^1,e^2,e^3$ (resp. $T^0,T^1,T^2,T^3$) in the present paper corresponds to $e_3,e_2,e_1,e_0$ in [ibid, p.105] (resp.\ $T^0,T^1,T_1,T_0$ in [ibid, Sect.\ 8.5.3]).
\qed
\medskip

{\it Proofs of (5A), (5B), (6A), (6B), (7A), and (7B) in 2.5.}

\medskip

We use mirror symmetry theorems (1)--(3) in Introduction.

From $e^{p}$ in (3B) and $\nabla_{\delta}$ in (4B), we produce $\tilde s^{p}$ inductively on $0\le p\le3$ as in (5B).
These are transported as (5A) in A-model.
For the last assertion in (5B) on the relation of $\nabla$-flat basis $\tilde s^p$ and the basis of the Frobenius solutions $y_j$, since $y_0T^3$, $y_0e^3=\Omega$, and the last expression in (5B) are killed by the operator $\cL$, it is enough to show the equality on the fiber $B_\bC:=\bC\otimes B$, i.e., the coincidence of the initial conditions.
We work in A-model.
By the asymptotics of the reverse relation of (5A) and of the expressions of $y_j$ in 2.2, we have
$$
\align
y_0T^3&=y_0\tilde s^3+y_0u\tilde s^2+y_0(2\pi i)^{-3}\frac{d\Phi}{du}\tilde s^1+y_0(2\pi i)^{-3}\Big(u\frac{d\Phi}{du}-2\Phi\Big)\tilde s^0\\
&\sim y_0\tilde s^3+y_0u\tilde s^2+\frac{5}{2}y_0u^2\tilde s^1
+\frac{5}{6}y_0u^3\tilde s^0\\
&\sim \tilde s^3+((2\pi i)^{-1}\log z)\tilde s^2+\frac{5}{2}((2\pi i)^{-1}\log z)^2\tilde s^1
+\frac{5}{6}((2\pi i)^{-1}\log z)^3\tilde s^0\\
&\sim y_0\tilde s^3+(2\pi i)^{-1}y_1\tilde s^2+5(2\pi i)^{-2}y_2\tilde s^1
+5(2\pi i)^{-3}y_3\tilde s^0.
\endalign
$$


To prove (6A), (6B), (7A), and (7B), we want to find $c^{10}, c^{21}, c^{20}, c^{32}, c^{31}, c^{30}\in\bC$ such that, on the fiber $B$,
$$
\align
\tilde s^0(u_0)&=b^0,\\
\tilde s^1(u_0)&=b^1+c^{10}\tilde s^0(u_0),\\
\tilde s^2(u_0)&=b^2+c^{21}\tilde s^1(u_0)+c^{20}\tilde s^0(u_0),\\
\tilde s^3(u_0)&=b^3+c^{32}\tilde s^2(u_0)+c^{31}\tilde s^1(u_0)+c^{30}\tilde s^0(u_0).
\endalign
$$
Then, since $\tilde s^{p}$ and $s^{p}$ are $\nabla$-flat, we have
$$
\align
s^0&=\tilde s^0,\\
s^1&=\tilde s^1-c^{10}\tilde s^0,\\
s^2&=\tilde s^2-c^{21}\tilde s^1-c^{20}\tilde s^0,\\
s^3&=\tilde s^3-c^{32}\tilde s^2-c^{31}\tilde s^1-c^{30}\tilde s^0.
\endalign
$$
Express the $e^{p}$ by the $s^{p}$ by using the inverse expressions of the above and of  (5B). 
Transporting these into A-model, we get expression of the $T^{p}$ by the $s^{p}$ and the $c^{jk}$.
Using Iritani asymptotics for the $s^{p}$ in A-model in 2.4, we get 
$$
\align
T^0&=\frac{1}{5}H^3,\\
T^1&=(2\pi i)^{-1}\frac{1}{5}H^2+\frac{c^{10}}{5}H^3,\\
T^2&\sim(2\pi i)^{-2}H+(2\pi i)^{-1}\Big(-\frac{1}{2}+\frac{c^{21}}{5}\Big)H^2\\
\quad&
+\Big(\Big(\frac{3}{2}-\frac{c^{21}}{5}+c^{10}\Big)u+\Big(\frac{7}{12}+\frac{c^{21}}{5}+\frac{c^{21}c^{10}}{5}+\frac{c^{20}}{5}\Big)\Big)H^3,\\
T^3&\sim(2\pi i)^{-3}+(2\pi i)^{-2}c^{32}H\\
\quad&
+(2\pi i)^{-1}\Big(\Big(-\frac{1}{2}+\frac{c^{21}}{5}-c^{32}\Big)u+\Big(\frac{5}{12}-\frac{c^{32}}{2}+\frac{c^{32}c^{21}}{5}+\frac{c^{31}}{5}\Big)\Big)H^2\\
\quad&
+\Big(\Big(\frac{1}{6}-\frac{c^{31}}{5}+\frac{c^{20}}{5}\Big)u+\Big(\frac{5i\zeta(3)}{\pi^3}+\frac{c^{30}}{5}\Big)\Big)H^3.
\endalign
$$
Since $T^p\in\cF^p$ by construction ((3A), (3B) in 2.5), 
we see that all coefficients of $H^j$ in the above 
expressions of $T^p$ are zero for $j+p>3$.
Thus we get
$$
\align
&c^{10}=-1, \;\; c^{21}=5/2,\;\; c^{20}=-35/12, \;\; \\
&c^{32}=0,\;\; c^{31}=-25/12,\;\;c^{30}=-25i\zeta(3)/\pi^3.
\endalign
$$
(6A), (6B), (7A), and (7B) follow from this and (5A), (5B).

We prove (8A), (8B).
The same argument goes for both cases.
We use the notation  in B-model. 
By the definition of the $\eta_j$ and (5B), 
$\Omega=\eta_0s^3+\eta_1s^2+\eta_2s^1+\eta_3s^0
=y_0\tilde s^3+(2\pi i)^{-1}y_1\tilde s^2+5(2\pi i)^{-2}y_2\tilde s^1+5(2\pi i)^{-3}y_3\tilde s^0$.
Substituting (6B) and comparing the coefficients of $\tilde s^p$, we have expressions of the $y_j$ by the $\eta_j$.
Solving these for the $\eta_j$, we get (8B).
\qed

\medskip

{\it Remark.}
It was pointed out by Hiroshi Iritani that the definitions and the descriptions of integral structures in [U14, 3.5, 3.6] are insufficient.
Actually, they were the first approximations of integral structures by means of $\gr^M$, 
which are characterized by the Frobenius solutions as in the last statement of (5A) and (5B) in 2.5.
The second proof in [ibid, 3.9] works well even in this approximation.

\bigskip

{\bf 3. Proof of (5) in Introduction}
\medskip

In this section, we prove (5.1)--(5.3) in Introduction for open mirror symmetry of quintic threefolds.
We prove them by constructing a normal function in log mixed Hodge theory for B-model in 3.1--3.2 below. 
This argument is applicable to the case of A-model by (1)--(4) in Introduction and the correspondence table in 2.5.
We give some discussions on geometries and local systems in 3.3.

\medskip

{\bf 3.1. Proofs of (5.1)--(5.3) in Introduction over log disc $S$}
\medskip

We consider B-model.
To make the monodromy of $\cT_B$ unipotent, we take a double cover $z^{1/2}\mapsto z$.
Let $S$ be a neighborhood disc of $p_0$ in the $z^{1/2}$-plane endowed with log 
structure associated to the divisor $p_0$ in $S$.
Denote by $\cH$ and $\cT$ the pullbacks of the log Hodge structure $\cH^{V^\circ}$ and the tension $\cT_B$ by the double covering, respectively.

We are looking for an extension $\tilde\cH$:
$$
0\to \cH\to \tilde\cH \to \bZ(-2)\to0
$$
of log mixed Hodge structures with $\gr^W_4\tilde\cH=\bZ(-2)$ and $\gr^W_3\tilde\cH=\cH$, which has 
liftings $1_\bZ$ and $1_F$ of $(2\pi i)^{-2}\cdot1\in\bZ(-2)$ in $\tilde\cH$ respecting the lattice and the Hodge filtration, respectively, such that the tension $\cT$ is described as
$$
Q(1_F-1_\bZ,\Omega)=\int_{C_-}^{C_+}\Omega=\cT, \tag1
$$
where $Q$ is the polarization of $\cH$ coming from 2.5 (1B) and $\Omega$ is the 3-form from 2.2.

To find such a log mixed Hodge structure, we use the basis $e^p$ $(0\le p\le 3)$ respecting the Deligne decomposition of $(\cH,\cM,\cF)$ from 2.5 (3B), and the $\nabla$-flat integral basis $s^p$ $(0\le p\le 3)$ from 2.5 (5B).
We also use the integral periods $\eta_j$ $(0\le j\le 3)$ in 2.5 (8B).
Note that these players  are already extended and live together over $S^{\log}$.

Let the local system $L_\bQ$ and the N\'eron model $J_{L_\bQ}$ be as in (5) in Introduction (see also 1.8).
Then $J_{L_\bQ}=\cE xt^1_{\LMH/S}(\bZ(-2),\cH)$, and let $1_\bZ:=((2\pi i)^{-2}\cdot1,-(\cT/\eta_0)s^0)\in \tilde\cH_\bZ$ be a lifting of $(2\pi i)^{-2}\cdot1\in\bZ(-2)=(\gr^W_4)_\bZ$, where $(\cT/\eta_0)s^0\in\cH_{\cO^{\log}}=(\gr^W_3)_{\cO^{\log}}$.
In particular, the connection $\nabla=\nabla^{\text{GM}}$ on $\cH$ is extended over $\tilde\cH$ by $\nabla(1_\bZ)=0$.

To find $1_F$, we write $1_F-1_\bZ=ae^3+be^2+ce^1-(\cT/\eta_0)e^0$ with
$a,b,c\in\cO^{\log}_S$ by using (1).
The Griffiths transversality condition on $1_F-1_\bZ$ is understood as vanishing 
of the coefficient of $e^0$ in $\nabla_\delta(1_F-1_\bZ)$.
Using 2.5 (4B), we have
$$
\nabla_\delta(1_F-1_\bZ)=
(\delta a)e^3+(a+\delta b)e^2
+\Big(b\frac{1}{(2\pi i)^3}\frac{d^3\Phi}{du^3}+\delta c\Big)e^1+(c-\delta(\cT/\eta_0))e^0.
$$
Hence, the above condition is equivalent to $c=\delta(\cT/\eta_0)$ and $a, b$ arbitrary.
Using the relation \lq\lq modulo $F^2$", we can take $a=b=0$.
Thus
$$
1_F=1_\bZ+(\delta(\cT/\eta_0))e^1-(\cT/\eta_0)e^0.
$$
The pair $1_\bZ$ and $1_F$ yields the desired element of $\cE xt^1_{\text{LMH}/S}(\bZ(-2),\cH)$, hence 
$1_F-1_\bZ$ yields the desired log normal function.
(5.1) is proved.

Next, we will find a splitting of the weight filtration $W$ of the local system $L_\bQ$, i.e., a splitting of $W$ which is compatible with the local monodromy of the local system $L_\bQ$.
We use the monodromy table in [W07, (3.14)].
This is computed for A-model but applicable also for B-model by the results (1)--(4) in Introduction and in Section 2.5.
Let $T_\infty^2$ be the monodromy around $p_0: z^{1/2}=0$ and put 
$N:=\log(T_\infty^2)$.
By [ibid], $N(\cT/\eta_0)=-1$ hence $N(1_\bZ)=s^0$.
On the other hand, we have $N(s^1)=-2s^0$.
(Here we use the rotation of the monodromy as $\log z \mapsto \log z + 2\pi i$.)
Define
$$
1_\bZ^{\text{spl}}:=1_\bZ+\frac{1}{2}s^1=\Big((2\pi i)^{-2}\cdot1,\,\frac{1}{2}s^1-\frac{\cT}{\eta_0}s^0\Big) \in \tilde\cH_\bQ.
\tag2
$$
Then $N(1_{\bZ}^{\spl})=0$, and this gives the desired splitting of $W$ of the local system $L_\bQ$.

A lifting $1_F^{\text{spl}}$ for $1_\bZ^{\text{spl}}$, respecting the Hodge filtration, is computed as before and we get
$$
1_F^{\text{spl}}=1_\bZ^{\text{spl}}+\Big(\delta\Big(\frac{\cT}{\eta_0}\Big)\Big)e^1-\frac{\cT}{\eta_0}e^0.
$$
The pair $1_\bZ^{\text{spl}}$ and $1_F^{\text{spl}}$ yields the desired element of $\cE xt^1_{\text{LMH}/S}(\bZ(-2),\cH)$ which splits the weight filtration $W$ of the local system $L_\bQ$.
Note that $1_F^{\text{spl}}-1_\bZ^{\text{spl}}=1_F-1_\bZ=(\delta(\cT/\eta_0))e^1-(\cT/\eta_0)e^0$.
(5.2) is proved.

(5.3) follows immediately from the above results.

We add a remark that the $W$-relative $N$-filtration $M=M(N,W)$ on $H_\bR$ in the admissibility condition 1.5 (1) is given by 
$$
\align
&M_{-1}=0\subset M_0=M_1=\bR s^0\subset M_2=M_3=M_1+\bR s^1\\
\subset\; &M_4=M_5=M_3+\bR s^2+\bR 1_{\bZ} \subset M_6=\tilde\cH_\bR=M_5+\bR s^3.
\endalign
$$
\medskip

{\bf 3.2. Proof of (5.1) and (5.2) in Introduction over log point $p_0$}
\medskip

We still consider B-model.
We show here that (5.1) and (5.2) in Introduction have meanings just over 
the log point $p_0$ and that the computations in their proofs become simpler.

Recall that
$$
\cT=-\frac{\eta_1}{2}-\frac{\eta_0}{4}+a_0\tau
\qquad \Big(a_0:=\frac{15}{\pi^2},\; \tau: \text{tau function}\Big)\tag1
$$
from [W07].
We substitute $z^{1/2}=0$ to $\cT$ carefully as follows.
Recall $\eta_1=\eta_0u$ from 2.2 and $u=x+iy$ from 2.4.
Write $v:=x+i\infty$ and define
$$
\cT(0):=-\frac{v}{2}-\frac{1}{4}+a_0
\qquad
\text{in $\cO_{p_0}^{\log}=\bC[v]$}.
$$
We abuse the notation $e^p$ and $s^p$ also for their restrictions over the log point $p_0$, and so they live together over $p_0^{\log}=(\bS^1,\bC[v])$.

Similarly as in (5) in Introduction, but using now $\cT(0)s^0$ instead of $(\cT/\eta^0)s^0$ because $\eta^0(0)=1$, we define a local system $L_\bQ$ and a N\'eron model $J_{L_\bQ}$ lying on $L_\bQ$.
Let $\tilde\cH$ be an extension of log mixed Hodge structures over the log point $p_0$, we are looking for, like in 3.1, and let $1_\bZ:=((2\pi i)^{-2}\cdot1,\,-\cT(0)s^0)$ be a lifting of $(2\pi i)^{-2}\cdot1\in \bZ(-2)=(\gr^W_4)_\bZ$ in $\tilde\cH_\bZ$. 
Hence the connection $\nabla$ on $\cH$ is extended over $\tilde\cH$ by $\nabla(1_\bZ)=0$.
Note that both $2\pi iq\frac{d}{dq}$ and $2\pi iz\frac{d}{dz}$ coincide with $\frac{d}{dv}$ 
now, which is denoted by $\delta$.
To find $1_F$, write $1_F-1_\bZ=ae_0+be_1+ce^1-\cT(0)e^0$ $(a,b,c\in\bC[v],\, \eta_0(0)=1)$ and compute $\nabla_\delta(1_F-1_\bZ)$ as in 3.1.
Then, by the Griffiths transversality, we have $c=-1/2$, $a$ and $b$ arbitrary.
By the relation \lq\lq modulo $F^2$", $a$ and $b$ can be reduced to $0$.
Thus, we have
$$
1_F=1_\bZ+(\delta\cT(0))e^1-\cT(0)e^0
=1_\bZ-\frac{1}{2}e^1+\Big(\frac{v}{2}+\frac{1}{4}-a_0\Big)e^0.
$$
The pair $1_\bZ$ and $1_F$ yields the desired element of $\cE xt^1_{\text{LMH}/S}(\bZ(-2),\cH)$.
(5.1) is proved.

The splitting of the weight filtration $W$ of the local system $L_\bQ$ is computed as 
in 3.1 but more simply, and we define
$$
1_\bZ^{\spl}:=1_\bZ+\frac{1}{2}s^1=\Big((2\pi i)^{-2}\cdot1,\,\frac{1}{2}s^1-\cT(0)s^0\Big) \in \tilde\cH_\bQ.
$$
Similarly, a lifting $1_F^{\spl}$ for $1_\bZ^{\spl}$ is computed simply, 
and we get
$$
1_F^{\spl}=1_\bZ^{\log}+(\delta\cT(0))e^1-\cT(0)e^0.
$$
(5.2) is proved.
\qed

\medskip

{\it Remark}.
Note that (5.3) does not have meaning in the present context.
This is because tau function disappears except its constant term when 
$z^{1/2}=0$ is substituted.
That is, in this step, we lose the transcendental data of the tension $\cT$, contained 
as the extension of its underlying local system, from which we can recover 
the position of the quintic mirror in its complex moduli space.

\bigskip


{\bf 3.3. Discussions on geometries and local systems}
\medskip

We discuss here the relation with geometries and local systems considered in [W07] and [MW09].
Forgetting Hodge structures, we consider only local systems corresponding to the monodromy of integral periods and tensions.

Let $V_\psi$ and $V^\circ_\psi$ be a quintic threefold and its mirror from 2.1.
Let $S$ be a small neighborhood in the $z$-plane ($z$ in 2.2) of the maximal unipotent monodromy point $p_0$ endowed with the log structure associated to the divisor $p_0$.

We first consider B-model.
Let the setting be as in [MW09, 4].
For $z\ne0$ near $0$, i.e., near $p_0$, let $V^{\circ}_z$ be the mirror quintic and $C_{+,z}\cup C_{-,z}$ be the disjoint union of smooth rational curves on $V^{\circ}_z$ coming from the two conics contained in $V_\psi\cap \{x_1+x_2=x_3+x_4=0\}\subset\bP^4(\bC)$.
From the relative homology sequence for $(V^{\circ}_z, (C_{+,z}\cup C_{-,z}))$, we have
$$
0@>>>H_3(V^\circ_z;\bZ)@>>>H_3(V^\circ_z, (C_{+,z}\cup C_{-,z});\bZ)@>\partial>>\bZ([C_{+,z}]-[C_{-,z}])@>>>0, \tag1
$$
where $\bZ([C_{+,z}]-[C_{-,z}])$ is $\Ker(H_2(C_{+,z}\cup C_{-,z});\bZ)\to H_2(V^{\circ}_z;\bZ))$.
The monodromy $T_\infty$ around $p_0$ interchanges $C_{+,z}$ and $C_{-,z}$.

Respecting the sequence (1), we take a family of cycles Poincar\'e duality isomorphic  to the flat integral basis $s^p$ ($0\le p\le 3$) in 2.4 and a family of chains joining from $C_{-,z}$ to $C_{+,z}$ (a choice up to integral cycles and up to half twists), and over them integrate the family of 3-forms $\Omega(z)$ with log pole over $z=0$ ($z$ in the punctured disc in the $z$-plan) in 2.2, then we have a family of vectors $(\eta_0,\eta_1,\eta_2,\eta_3,\cT)$ consisting of periods (2.5 (8B)) and a tension.
This corresponds to the data in [W07], [MW09] (cf.\ Remark in 2.5).
Since $T_\infty(\cT)=-(\cT+\eta_1+\eta_0)$ by [W07, (3.14)], we find $\cT+\frac{1}{2}\eta_1+\frac{1}{4}\eta_0=\frac{15}{\pi^2}\tau$ (see 3.2 (1)) is an eigenvector of the monodromy $T_\infty$ with eigenvalue $-1$.

The family of sequences (1) ($z\ne0$) forms an exact sequence of local systems of $\bZ$-modules.
Pulling this back to $S^*$ in 3.1 by the double cover $z^{1/2}\mapsto z$, we have a sequence with unipotent local monodromy and its extension over $S^{\log}$.
Applying Tate twist $(-3)$ and Poincar\'e duality isomorphism to the left and the right ends of this exact sequence, we have a local system $L'$ over $S^{\log}$ which is an extension of $\bZ(-2)$ by $\cH_\bZ$:
$$
0@>>>\cH_\bZ@>>>L'@>>>\bZ(-2)@>>>0.\tag2
$$
Take a lifting $1_\bZ:=((2\pi i)^{-2}\cdot1,\,(\cT/\eta_0)s^0)$ in $L'$ of $(2\pi i)^{-2}\cdot1\in\bZ(-2)$, and extend $\nabla$ on $\cH_\bZ$ over $L'$ by $\nabla(1_\bZ)=0$.
We look for a $\nabla$-flat $T_\infty^2$-invariant element associated to $1_\bZ$.
This is done as in 3.1, and we get $1_\bZ^{\text{spl}}:=1_\bZ-(s^1/2)$.
Thus we know that $L'$ coincides with $\tilde\cH_\bZ$ in 3.1.

For the relative monodromy weight filtration $M=M(N,W)$ of $L^\prime$, we see that $1_\bZ\in M_4$ and $s^1\in M_2$ are the smallest filters containing each element in question.
Taking the graded quotients by $M$ of the sequence (2), we have
$$
\align
&\gr^M_6\cH_\bZ@>\sim>>\gr^M_6L',\tag3\\
0@>>>&\gr^M_4\cH_\bZ@>>>\gr^M_4L'@>>>\bZ(-2)@>>>0,\\
0@>>>&\gr^M_2\cH_\bZ@>>>\gr^M_2L'@>>>(\text{$2$-torsion})@>>>0,\\
&\gr^M_0\cH_\bZ@>\sim>>\gr^M_0L'.
\endalign
$$
Here we abuse the notation $M$ also for the monodromy filtration on $\cH_\bZ$, because it coincides with the restriction of $M=M(N,W)$ to $\cH_\bZ$.
The 2-torsion in the third sequence of (3) corresponds to a half twist of chains from $C_-$ to $C_+$.
Standing on a half integral point and looking at the integral points nearby, we have two orientations.
These correspond to the two orientations of a half twist of the chains, and also correspond to $\cT_{\pm}:=\pm(\frac{15}{\pi^2}\tau-\frac{\eta_0}{4})-\frac{\eta_1}{2}$ in [W07].
$\cT_-$ is different from $-\cT_+$ by the complementary half twist, i.e., $\cT_++\cT_-=-\eta_1$.

By using mirror symmetries (1)--(5), or more precisely, by the results in Section 2.5 and Section 3.1, 
$\cH_\bZ=\cH_\bZ^{V^\circ}$, $\cT=\cT_B$, $1_\bZ$, $\nabla=\nabla^{\text{GM}}$, $1_\bZ^{\spl}$, and $M=M(W,N)$ of B-model are transformed to 
the corresponding $\cH_\bZ=\cH_\bZ^{V}$, $\cT=y_0\cT_A$, $1_\bZ$, $\nabla=\nabla^{\text{even}}$, $1_\bZ^{\spl}$, and $M=M(W,N)$ of A-model, and the exact sequences (2) and (3) of B-model are transformed to the corresponding exact sequences of A-model.

It is interesting to study the relations of these exact sequences with the geometries of Fermat quintic $V=V_\psi$ with $\psi=0$ and its Lagrangian submanifold $Lg:=V\cap\bP^4(\bR)$ in [W07, 2.1], [MW09, 3].

\medskip

{\it Remark.}
(1) The argument in 3.3 can be performed even over the log point $p_0$.

(2) [PSW08] and [LiLY12] are related with the topics in this subsection.



\vskip20pt




\Refs

\widestnumber\key{COGP91}


\ref
\key CoK99
\by  D.\ A.\ Cox and S.\ Katz
\book Mirror symmetry and algebraic geometry
\bookinfo Math\.Surveys and Monographs
\vol 68
\yr 1999
\pages 469. MR 2000d:14048
\publaddr AMS
\endref

\ref
\key CDGP91
\by P\.Candelas, C\.de la Ossa, P\.S\.Green, and L\.Parks
\paper A pair of Calabi-Yau manifolds as an exactly soluble 
superconformal theory
\jour Nuclear Physics
\vol B 358
\yr 1991
\pages 21-74
\endref

\ref
\key D70
\by P. Deligne
\book \'Equations differentielles $\grave{a}$ points singuliers
r\'eguliers
\publ Lect. Notes in Math. No. 163, Springer-Verlag
\yr 1970. MR 54$\sharp$5232
\endref

\ref
\key D97
\bysame
\paper Local behavior of Hodge structures at infinity
\jour in Mirror Symmetry II (B.\ Greene and S.-T.\ Yau, eds.), 
AMS/IP Stud.\ Adv.\ Math.\ {\bf1}, 1997
\vol 
\yr 
\pages 683--699. MR 98a:14015
\endref

\ref
\key G96
\by A.\ B.\ Givental
\paper Equivariant Gromov--Witten invariants
\jour Internat.\ Math.\ Res.\ Notes
\vol 13
\yr 1996
\pages 613--663
\endref

\ref
\key I09
\by H.\ Iritani
\paper An integral structure in quantum cohomology
\jour Adv.\ Math.\ 
\vol 222 (3)
\yr 2009
\pages 1016--1079
\endref

\ref
\key I11
\bysame
\paper Quantum cohomology and periods
\jour Ann.\ Inst.\ Fourier (Grenoble)
\vol 61 no.7
\yr 2011
\pages 2909 -- 2958 
\endref

\ref
\key KN99
\by  K\.Kato and C.\ Nakayama
\paper Log Betti cohomology, log \'etale cohomology, and log de Rham cohomology of log schemes over $\bC$
\jour Kodai Math. J.
\vol 22
\yr 1999
\pages 161--186. MR 2000i:14023
\endref

\ref
\key KNU08
\by K.\ Kato, C.\ Nakayama and S.\ Usui
\paper $\SL(2)$-orbit theorem for degeneration of mixed Hodge structure
\jour J.\ Algebraic Geometry
\vol 17
\yr 2008
\pages 401--479. MR 2009b:14020
\endref

\ref
\key KNU09
\bysame
\paper Classifying spaces of degenerating mixed Hodge structures, I\rom:
Borel--Serre spaces
\jour Advanced Studies in Pure Math. {\bf54}:
Algebraic Analysis and Around, 2009
\yr 
\pages 187--222. MR 2010g:14010
\endref

\ref
\key KNU11
\bysame
\paper Classifying spaces of degenerating mixed Hodge structures, II\rom:
Spaces of $\SL(2)$-orbits
\jour Kyoto J.\ Math.\ {\bf 51-1}: Nagata Memorial Issue
\yr 2011
\pages 149--261. MR 2012f:14012
\endref

\ref
\key KNU13
\bysame
\paper Classifying spaces of degenerating mixed Hodge structures, III\rom: Spaces of nilpotent orbits
\jour J.\ Algebraic Geometry, 
\vol 22
\yr 2013
\pages 671--772
\endref

\ref
\key KNU14
\bysame
\paper N\'eron models for admissible normal functions
\jour Proc.\ Japan Academy
\vol 90{\rm, Ser.\ A}
\yr 2014
\pages 6--10
\endref


\ref
\key KU99
\by K.\ Kato and S.\ Usui
\paper Logarithmic Hodge structures and classifying
spaces {\rm(summary)}
\jour in CRM Proc.\ \& Lect.\ Notes:
The Arithmetic and Geometry of Algebraic Cycles,
(NATO Advanced Study Institute /
CRM Summer School 1998: Banff, Canada)
\vol 24
\yr 1999
\pages 115--130. MR 2001e:14009
\endref

\ref
\key KU02
\bysame
\paper Borel-Serre spaces and spaces of
{\rm SL(2)}-orbits
\jour Advanced Studies in Pure Math.\ {\bf36}:
Algebraic Geometry 2000, Azumino,
\yr 2002
\pages 321--382. MR 2004f:14021
\endref

\ref
\key KU09
\bysame
\book Classifying spaces of degenerating polarized 
Hodge structures
\bookinfo Ann\. Math\. Studies, Princeton Univ\.Press
\vol 169
\yr 2009
\pages 288. MR 2009m:14012
\publaddr Princeton
\endref

\ref
\key LW12
\by G.\ Laporte and J.\ Walcher
\paper  Monodromy of an inhomogenoues Picard--Fuchs equation
\jour SIGMA 
\vol 8
\yr 2012
\pages 056, 10 pages
\endref

\ref
\key LiLY12
\by S.\ Li, B.\ Lian, and S.-T.\ Yau
\paper Picard--Fuchs equations for relative periods and Abel-Jacobi map for Calabi-Yau hypersurfaces
\jour Amer.\ J.\ Math.\
\vol 134-5
\yr 2012
\pages 1345-1384. MR 2975239
\endref

\ref
\key LLuY97
\by B.\ Lian, K.\ Liu, and S.-T.\ Yau
\paper Mirror principle I
\jour Asian J.\ Math.\
\vol 1
\yr 1997
\pages 729--763. MR 99e:14062
\endref

\ref
\key M93
\by D.\ Morrison
\paper Mirror symmetry and rational curves on quintic threefolds:
A guide for mathematicians
\jour J. of AMS
\vol 6-1
\yr 1993
\pages 223--247. MR 93j:14047
\endref

\ref
\key M97
\bysame
\paper Mathematical aspects of mirror symmetry
\jour in Complex algebraic geometry (Park City, UT, 1993), IAS/Park City Math.\ Ser.\ {\bf 3}, AMS
\vol 
\yr 1997
\pages 265--327. MR 98g:14044
\endref

\ref
\key MW09
\by D.\ Morrison and J.\ Walcher
\paper D-branes and normal functions
\jour Adv.\ Theor.\ Math.\ Phys.\
\vol 13-2
\yr 2009
\pages 553--598. MR 2010b:14081
\endref

\ref
\key O03
\by A. Ogus
\paper On the logarithmic Riemann-Hilbert correspondences
\jour Documenta Math. Extra volume:
Kazuya Kato's Fiftieth birthday
\yr 2003
\pages 655--724
\endref

\ref
\key P98
\by R.\ Pandharipande
\paper Rational curves on hypersurfaces {\rm[after A.\ Givental]}
\jour S\'eminaire Bourbaki 848, Ast\'erisque
\vol 252
\yr 1998
\pages Exp.\ No.\ 848, 5, 307--340. MR 2000e:14094
\endref

\ref
\key PSW08
\by R.\ Pandharipande, J.\ Solomon, and J.\ Walcher
\paper Disk enumeration on the quintic 3-fold
\jour J.\ Amer.\ Math.\ Soc.\
\vol 21-4
\yr 2008
\pages 1169--1209. MR 2009j:14075
\endref

\ref
\key S73
\by W.\ Schmid
\paper Variation of Hodge structure\rom:
The singularities of the period mapping
\jour Invent.\ Math.
\vol 22
\yr 1973
\pages 211--319. MR 52$\sharp$3157
\endref

\ref
\key U84
\by S.\ Usui
\paper Variation of mixed Hodge structure arising from
family of logarithmic deformations II\rom: Classifying space
\jour Duke Math\.J.
\vol 51-4
\yr 1984
\pages 851--875. MR 86h:14005
\endref

\ref
\key U08
\bysame
\paper Generic Torelli theorem for quintic-mirror family
\jour Proc.\ Japan Acad.
\vol 84, Ser.\ A, No.\ 8
\yr 2008
\pages 143--146. MR 2010b:14012
\endref

\ref
\key U14
\bysame
\paper A study of mirror symmetry through log mixed Hodge theory
\jour Hodge Theory, Complex Geometry, and Representation Theory, Contemporary Math., AMS 
\vol  608
\yr 2014
\pages 285--311
\endref

\ref
\key W07
\by J.\ Walcher
\paper Opening mirror symmetry on the quintic
\jour Commun.\ Math.\ Phys.\
\vol 276
\yr 2007
\pages 671--689. MR 2008m:14111
\endref





\endRefs
\enddocument